\documentclass{siamart171218}

\usepackage{amsfonts}
\usepackage{mathrsfs}
\usepackage{mathtools}
\usepackage{graphicx}
\usepackage{epstopdf}
\usepackage{algorithmic}
\usepackage{comment}
\usepackage{xcolor}
\usepackage{tabto}
\usepackage{multirow}
\usepackage[caption=false]{subfig}
\usepackage{bm}
\usepackage{tikz}

\newsiamremark{remark}{Remark}

\crefname{algocf}{Algorithm}{Algorithms}
\Crefname{algocf}{Algorithm}{Algorithms}
\crefname{subsection}{section}{sections}
\Crefname{subsection}{Section}{Sections}

\renewcommand{\theequation}{\thesection.\arabic{equation}}
\renewcommand{\thefigure}{\thesection.\arabic{figure}}
\renewcommand{\thetable}{\thesection.\arabic{table}}

\numberwithin{theorem}{section}
\numberwithin{equation}{section}
\numberwithin{figure}{section}
\numberwithin{table}{section}
\numberwithin{algorithm}{section}

\title{A Low-rank Solver for the Stochastic Unsteady Navier--Stokes Problem%
   \thanks{This work was supported by the U.S. Department of Energy Office of Advanced Scientific Computing Research, Applied Mathematics program under award DE-SC0009301 and by the U.S. National Science Foundation under grant DMS1819115.}}

\author{
   Howard C. Elman%
   \thanks{Department of Computer Science and Institute for Advanced Computer Studies, University of Maryland, College Park, MD 20742 (\email{elman@cs.umd.edu}).}
   \and
   Tengfei Su%
   \thanks{Applied Mathematics \& Statistics, and Scientific Computation Program, University of Maryland, College Park, MD 20742 (\email{tengfesu@math.umd.edu}).}
}

\headers{Low-rank Solver for Stochastic Unsteady Navier--Stokes Problem}
{H. C. Elman and T. Su}

\begin{document}

\maketitle

\begin{abstract}
We study a low-rank iterative solver for the unsteady Navier--Stokes equations for incompressible flows with a stochastic viscosity. The equations are discretized using the stochastic Galerkin method, and we consider an all-at-once formulation where the algebraic systems at all the time steps are collected and solved simultaneously. The problem is linearized with Picard's method. To efficiently solve the linear systems at each step, we use low-rank tensor representations within the Krylov subspace method, which leads to significant reductions in storage requirements and computational costs. Combined with effective mean-based preconditioners and the idea of inexact solve, we show that only a small number of linear iterations are needed at each Picard step. The proposed algorithm is tested with a model of flow in a two-dimensional symmetric step domain with different settings to demonstrate the computational efficiency.
\end{abstract}


%-------------------------------------------------------------
\section{Introduction}
\label{sec:ns}

Stochastic partial differential equations (PDEs) are widely used to model physical problems with uncertainty \cite{LeKn10}. In this paper, we develop some new computational methods for solving the stochastic unsteady Navier--Stokes equations, using stochastic Galerkin methods \cite{GhSp03} to address the stochastic nature of the problem and so-called all-at-once treatment of time integration.Ê

For a time-dependent problem, the solutions at different time steps are usually computed in a sequential manner via time stepping. For example, a fully-implicit scheme with adaptive time step sizes was studied in \cite{ElMi11,KaGr10}. On the other hand, an all-at-once system can be formed by collecting the algebraic systems at all the discrete time steps into a single one, and the solutions are computed simultaneously. Such a formulation avoids the serial nature of time stepping, and allows parallelization in the time direction for accelerating the solution procedure \cite{GaNe16,MaRo08,McPe18}. A drawback, however, is that for large-size problems, the all-at-once system may require excessive storage. In this study, we address this issue by using a low-rank tensor representation of data within the solution methods.

We develop a low-rank iterative algorithm for solving the unsteady Navier--Stokes equations with an uncertain viscosity. The equations are linearized with Picard's method. At each step of the nonlinear iteration, the stochastic Galerkin discretization gives rise to a large linear system, which is solved by a Krylov subspace method. Similar approaches have been used to study the steady-state problem \cite{PoSi12, SoEl16b}, where the authors also proposed effective preconditioners by taking advantage of the special structures of the linear systems. To reduce memory and computational costs, we compute low-rank approximations to the discrete solutions, which are represented as three-dimensional tensors in the all-at-once formulation. We refer to \cite{GrKr13} for a review of low-rank tensor approximation techniques, and we will use the tensor train decomposition \cite{Os11} in this work. The tensor train decomposition allows efficient basic operations on tensors. A truncation procedure is also available to compress low-rank tensors in the tensor train format to ones with smaller ranks. 

Our goal is to use the low-rank tensors within Krylov subspace methods, in order to efficiently solve the large linear systems arising in each nonlinear step. The basic idea is to represent all the vector quantities that arise during the course of a Krylov subspace computation as low-rank tensors. With this strategy, much less memory is needed to store the data produced during the iteration. Moreover, the associated computations, such as matrix-vector products and vector additions, become much cheaper. The tensors are compressed in each iteration to maintain low ranks. This idea has been used for the conjugate gradient (CG) method and the generalized minimal residual (GMRES) method, with different low-rank tensor formats \cite{AnTo15, BaGr13, Do13, KrTo11}. In addition, the convergence of Krylov subspace methods can be greatly improved by an effective preconditioner. In conjunction with the savings achieved through low-rank tensor computations, we will derive preconditioners for the stochastic all-at-once formulation based on some state-of-the-art techniques used for deterministic problems, and we will demonstrate their performances in numerical experiments. We also explore the idea of inexact Picard methods where the linear systems are solved inexactly at each Picard step to further save computational work, and we show that with this strategy very small numbers of iterations are needed for the Krylov subspace method.

We note that a different type of approach, the alternating iterative methods \cite{DoSa14,HoRo12,Sc05}, including the density matrix renormalization group (DMRG) algorithm and its variants, can be used for solving linear systems in the tensor train format. In these methods, each component of the low-rank solution tensor is approached directly and optimized by projecting to a small local problem. This approach avoids the rank growth in intermediate iterates typically encountered in a low-rank Krylov subspace method. However, these methods are developed for solving symmetric positive definite systems and require nontrivial effort to be adapted for a nonsymmetric Navier--Stokes problem.

The rest of the paper is organized as follows. In \cref{sec:setting} we give a formal presentation of the problem. Discretization techniques that result in an all-at-once linear system at each Picard step are discussed in \cref{sec:discrete}. In \cref{sec:low_rank} we introduce the low-rank tensor approximation and propose a low-rank Krylov subspace iterative solver for the all-at-once systems. The preconditioners are derived in \cref{sec:prec} and numerical results are given in \cref{sec:num_test}.

%-------------------------------------------------------------
\section{Problem setting}
\label{sec:setting}

Consider the unsteady Navier--Stokes equations for incompressible flows on a space-time domain $\mathcal{D}\times (0,t_f]$,
\begin{equation}
\label{eq:ns_0}
\begin{aligned}
    \frac{\partial \vec{u}}{\partial t} - \nabla\cdot(\nu\nabla\vec{u}) + \vec{u}\cdot\nabla\vec{u} + \nabla p &= \vec{0}, \\
    \nabla\cdot\vec{u} &= 0,
\end{aligned}
\end{equation}
where $\vec{u}$ and $p$ stand for the velocity and pressure, respectively, $\nu$ is the viscosity, and $\mathcal{D}$ is a two-dimensional spatial domain with boundary $\partial\mathcal{D}=\partial\mathcal{D}_\text{D}\cup\partial\mathcal{D}_\text{N}$. The Dirichlet boundary $\partial\mathcal{D}_\text{D}$ consists of an inflow boundary and fixed walls, and Neumann boundary conditions are set for the outflow,
\begin{equation}
\begin{aligned}
    \vec{u} &= \vec{u}_\text{D} & \text{ on }\partial{\mathcal{D}}_\text{D},\\
    \nu\nabla\vec{u}\cdot\vec{n} - p\vec{n} &= \vec{0} & \text{ on }\partial{\mathcal{D}}_\text{N}.
\end{aligned}
\end{equation}
We assume the Neumann boundary $\partial\mathcal{D}_\text{N}$ is not empty so that the pressure $p$ is uniquely determined. The function $\vec{u}_\text{D}(x,t)$ denotes a time-dependent inflow, typically growing from zero to a steady state, and it is set to zero at fixed walls. The initial conditions are zero everywhere for both $\vec{u}$ and $p$. 

The uncertainty in the problem is introduced by a stochastic viscosity $\nu$, which is modeled as a random field depending on a finite collection of random variables $\{\xi_l\}_{l=1}^{m}$ (or written as a vector $\xi$). Specifically, we consider a representation as a truncated Karhunen--Lo\`{e}ve (KL, \cite{Loeve}) expansion,
\begin{equation}
    \nu(x,\xi) = \nu_0(x) + \sum_{l=1}^{m} \nu_l(x)\xi_l,
\end{equation}
where $\nu_0$ is the mean viscosity, and $\{\nu_l\}_{l=1}^{m}$ are determined by the covariance function of $\nu$. We assume that the random parameters $\{\xi_l\}_{l=1}^{m}$ are independent and that the viscosity satisfies $\nu(x,\xi)\geq\nu_\text{min}>0$ almost surely for any $x\in\mathcal{D}$. We refer to \cite{PoSi12, SoEl16b} for different forms of the stochastic viscosity. The solutions $\vec{u}$ and $p$ in \cref{eq:ns_0} will also be random fields which depend on the space parameter $x$, time $t$, and the random variables $\xi$.

%-------------------------------------------------------------
\section{Discrete problem}
\label{sec:discrete}

In this section, we derive a fully discrete problem for the stochastic unsteady Navier--Stokes equations \cref{eq:ns_0}. This involves a time discretization scheme and a stochastic Galerkin discretization for the physical and parameter spaces at each time step. The discretizations give rise to a nonlinear algebraic system. Instead of solving such a system at each time step, we collect the systems from all time steps to form an all-at-once system, where the discrete solutions at all the time steps are solved simultaneously. The discrete problem is then linearized with Picard's method, and a large linear system is solved at each step of the nonlinear iteration.

\subsection{Time discretization}
For simplicity we use the backward Euler method for time discretization, which is first-order accurate but unconditionally stable and dissipative. The all-at-once formulation discussed later in \cref{sec:all_at_once} requires predetermined time steps. Divide the interval $(0,t_f]$ into $n_t$ uniform steps $\{t_k\}_{k=1}^{n_t}$ with step size $\tau=t_f/n_t$ and initial time $t_0=0$. Given the solution at time $t_{k-1}$, we need to solve the following equations for $\vec{u}^k$ and $p^k$:
\begin{equation}
\label{eq:ns_dt}
\begin{aligned}
    \frac{\vec{u}^{k}-\vec{u}^{k-1}}{\tau} - \nabla\cdot(\nu\nabla\vec{u}^k) + \vec{u}^k\cdot\nabla\vec{u}^k + \nabla p^k &= \vec{0}, \\
    \nabla\cdot\vec{u}^k &= 0.
\end{aligned}
\end{equation}
After discretization (in physical space and parameter space) the implicit method requires solving an algebraic system at each time step. In the following we discuss how the system is assembled from the stochastic Galerkin discretization of \cref{eq:ns_dt}.

\subsection{Stochastic Galerkin method}
At time step $k$, the stochastic Galerkin method finds parametrized approximate velocity solutions $\vec{u}_h^k$ and pressure solutions $p_h^k$ in finite-dimensional subspaces of $(H^1(\mathcal{D}))^2\otimes L^2(\Gamma)$ and $L^2(\mathcal{D})\otimes L^2(\Gamma)$, where $\Gamma$ is the joint image of the random variables $\{\xi_l\}$. The functional spaces are defined as follows,
\begin{equation}
\begin{aligned}
    (H^1(\mathcal{D}))^2\otimes L^2(\Gamma) &\coloneqq \Big\{ \vec{v}:\mathcal{D}\times\Gamma\rightarrow\mathbb{R} \mid \mathbb{E}\big[ \|\vec{v}\|_{(H^1(\mathcal{D}))^2}^2 \big] < \infty \Big\}, \\
    L^2(\mathcal{D})\otimes L^2(\Gamma) &\coloneqq \Big\{ q:\mathcal{D}\times\Gamma\rightarrow\mathbb{R} \mid \mathbb{E}\big[ \|q\|_{L^2(\mathcal{D})}^2 \big] < \infty \Big\}. \\
\end{aligned}
\end{equation}
The expectations are taken with respect to the joint distribution of the random variables $\{\xi_l\}$. In the following we use $\langle\cdot\rangle$ to denote the expected value. Let the finite-dimensional subspaces be $\mathcal{X} = \text{span}\{\vec{\phi}_i(x)\} \subset (H^1(\mathcal{D}))^2$, $\mathcal{Y} = \text{span}\{\varphi_i(x)\} \subset L^2(\mathcal{D})$, and $\mathcal{Z} = \text{span}\{\psi_r(\xi)\} \subset L^2(\Gamma)$. Let $\mathcal{X}_\text{D}^k$ and $\mathcal{X}_0$ be the spaces of functions in $\mathcal{X}$ with Dirichlet boundary conditions $\vec{u}_\text{D}(x,t_k)$ and $\vec{0}$ imposed for the velocity field, respectively. Then for \cref{eq:ns_dt} the stochastic Galerkin formulation entails the computation of $\vec{u}_h^k\in\mathcal{X}_\text{D}^k\otimes\mathcal{Z}$ and $p_h^k\in\mathcal{Y}\otimes\mathcal{Z}$, satisfying the weak form
\begin{equation}
\begin{aligned}
    \tau^{-1}\langle(\vec{u}_h^k,\vec{v}_h)\rangle - \tau^{-1}\langle(\vec{u}_h^{k-1},\vec{v}_h)\rangle + \langle(\nu\nabla\vec{u}_h^k,\nabla\vec{v}_h)\rangle \hspace{2.4cm}& \\ +\, \langle(\vec{u}_h^k\cdot\nabla\vec{u}_h^k,\vec{v}_h)\rangle - \langle(p_h^k,\nabla\cdot\vec{v}_h)\rangle &= 0, \\
    \langle(\nabla\cdot\vec{u}_h^k,q_h)\rangle &= 0,
\end{aligned}
\end{equation}
for any $\vec{v}_h\in\mathcal{X}_0\otimes\mathcal{Z}$ and $q_h\in\mathcal{Y}\otimes\mathcal{Z}$. Here, $(\cdot,\cdot)$ denotes the inner product in ${L}^2(\mathcal{D})$. For the physical spaces, we use a div-stable Taylor--Hood discretization \cite{ElSi14} on quadrilateral elements, with biquadratic basis functions $\{\vec{\phi}_i\}_{i=1}^{n_u} = \left\{ \left(\begin{smallmatrix}\phi_i\\0 \end{smallmatrix}\right), \left(\begin{smallmatrix}0\\ \phi_i \end{smallmatrix}\right) \right\}_{i=1}^{n_u/2}$ for velocity, and bilinear basis functions $\{\varphi_i\}_{i=1}^{n_p}$ for pressure. The stochastic basis functions $\{\psi_r\}_{r=1}^{n_\xi}$ are $m$-dimensional orthonormal polynomials constructed from generalized polynomial chaos (gPC, \cite{XiKa02}) satisfying $\langle \psi_r\psi_s \rangle=\delta_{rs}$. The stochastic Galerkin solutions are expressed as linear combinations of the basis functions,
\begin{equation}
\label{eq:sg_sol}
\begin{aligned}
    \vec{u}_{h}^k(x,\xi) &= \sum_{s=1}^{n_\xi}\sum_{j=1}^{n_u} u_{js}^k \vec{\phi}_j(x) \psi_s(\xi), \\
    p_{h}^k(x,\xi) &= \sum_{s=1}^{n_\xi}\sum_{j=1}^{n_p} p_{js}^k \varphi_j(x) \psi_s(\xi).
\end{aligned}
\end{equation}
The coefficient vectors $\bm{u}^k=[u_{11}^k,u_{21}^k,\ldots,u_{n_u1}^k,\ldots,u_{1n_\xi}^k,u_{2n_\xi}^k,\ldots,u_{n_un_\xi}^k]$ and similarly defined $\bm{p}^k$ are computed from the nonlinear algebraic system
\begin{equation}
\label{eq:ns_sys_k}
    \begin{pmatrix}
    \mathbb{F}^k(\bm{u}) & I_{n_\xi}\otimes B^T\\
    I_{n_\xi}\otimes B & 0
    \end{pmatrix}
    \begin{pmatrix} \bm{u}^k \\ \bm{p}^k \end{pmatrix} +
    \begin{pmatrix}
    -\tau^{-1}(I_{n_\xi}\otimes\bm{M}) & 0\\
    0 & 0
    \end{pmatrix}
    \begin{pmatrix} \bm{u}^{k-1}\\ \bm{p}^{k-1} \end{pmatrix} = 
    \begin{pmatrix} \bm{f}^{u,k}\\ \bm{f}^{p,k} \end{pmatrix}
\end{equation}
where
\begin{equation}
\label{eq:F^k}
    \mathbb{F}^k(\bm{u}) = \tau^{-1}(I_{n_\xi}\otimes\bm{M}) + \sum_{l=0}^{m}(G_l\otimes\bm{A}_l) + \sum_{l=1}^{n_\xi}(H_l\otimes\bm{N}(\vec{u}_{h,l}^k)).
\end{equation}
Here $I_{n_\xi}$ is the $n_\xi\times n_\xi$ identity matrix, and $\otimes$ denotes the Kronecker product of two matrices. The boldface matrices $\bm{M}$, $\bm{A}_l$, and $\bm{N}(\vec{u}_{h,l}^k)$ are $2\times 2$ block-diagonal, with the scalar mass matrix $M$, weighted stiffness matrix $A_l$, and discrete convection operator $N(\vec{u}_{h,l}^k)$ as diagonal components, where
\begin{equation}
\label{eq:mat_entries}
    [M]_{ij}=(\phi_j,\phi_i),\,\, [A_l]_{ij}=(\nu_l\nabla\phi_j,\nabla\phi_i),\,\, [N(\vec{u}_{h,l}^k)]_{ij}=(\vec{u}_{h,l}^k\cdot\nabla\phi_j,\phi_i),
\end{equation}
for $i,j=1,\ldots,n_u/2$. Note the dependency on $\bm{u}^k$ comes from the nonlinear convection term $\bm{N}$, with convection velocity $\vec{u}_{h,l}^k=\sum_j u_{jl}^k\vec{\phi}_j(x)$. Let $x=(x_1,x_2)$. The discrete divergence operator $B=[B_{x_1},B_{x_2}]$, with
\begin{equation}
    [B_{x_1}]_{ij} = -(\varphi_i,\frac{\partial\phi_j}{\partial x_1}),\,\,[B_{x_2}]_{ij} = -(\varphi_i,\frac{\partial\phi_j}{\partial x_2}),
\end{equation}
for $i=1,\ldots,n_p$ and $j=1,\ldots,n_u/2$. The matrices $\{G_l\}_{l=0}^{m}$ and $\{H_l\}_{l=1}^{n_\xi}$ of \cref{eq:F^k} come from the stochastic basis functions and have entries 
\begin{equation}
    [G_l]_{rs} = \langle \xi_l\psi_r\psi_s \rangle,\,\, [H_l]_{rs} = \langle \psi_l\psi_r\psi_s \rangle,
\end{equation}
for $r,s=1,\ldots,n_\xi$, where $\xi_0\equiv 1$. These matrices are also sparse due to orthogonality of the basis functions \cite{ErUl10}. The Dirichlet boundary conditions are incorporated in the right-hand side of \cref{eq:ns_sys_k}.

\subsection{All-at-once system}
\label{sec:all_at_once}
As discussed in the beginning of the section, we consider an all-at-once system where the discrete solutions at all the time steps are computed together. Let
\begin{equation}
    \bm{u} = \begin{pmatrix}
    \bm{u}^1 \\ 
    \bm{u}^2 \\
    \vdots \\
    \bm{u}^{n_t} \end{pmatrix}  \in \mathbb{R}^{n_tn_\xi n_u}
\end{equation}
and let $\bm{p}$, $\bm{f}^u$, and $\bm{f}^p$ be similarly defined. By collecting the algebraic systems \cref{eq:ns_sys_k} corresponding to all the time steps $\{t_k\}_{k=1}^{n_t}$, we get the single system
\begin{equation}
\label{eq:ns_sys_all}
    \begin{pmatrix}
    \mathbb{F}(\bm{u})+\mathbb{C} & \mathbb{B}^T\\
    \mathbb{B} & 0
    \end{pmatrix}
    \begin{pmatrix} \bm{u}\\ \bm{p} \end{pmatrix} =
    \begin{pmatrix} \bm{f}^u\\ \bm{f}^p \end{pmatrix},
\end{equation}
where $\mathbb{F}(\bm{u})$ is block diagonal with $\mathbb{F}^k(\bm{u})$ as the $k$th diagonal block, $\mathbb{B}=I_{n_t}\otimes I_{n_\xi}\otimes B$, and $\mathbb{C}=-\tau^{-1}C_{n_t}\otimes I_{n_\xi}\otimes\bm{M}$ with $C_{n_t} = \left( \begin{smallmatrix} 0&&&\\ 1&0&&\\ &\ddots&\ddots&\\ &&1&0 \end{smallmatrix} \right)\in\mathbb{R}^{n_t\times n_t}$. Note that the zero initial conditions are incorporated in \cref{eq:ns_sys_k} for $k=1$. The all-at-once system \cref{eq:ns_sys_all} is nonsymmetric and blockwise sparse. Each part of the system contains sums of Kronecker products of three matrices, i.e., in the form $\sum_l X_l^{(1)}\otimes X_l^{(2)}\otimes X_l^{(3)}$. In fact, from \cref{eq:F^k},
\begin{equation}
\label{eq:F}
    \mathbb{F}(\bm{u})=\tau^{-1}I_{n_t}\otimes I_{n_\xi}\otimes\bm{M} + \sum_{l=0}^{m}(I_{n_t}\otimes G_l\otimes\bm{A}_l)+\mathbb{N}(\bm{u}).
\end{equation}
We discuss later (see \cref{sec:stoch_conv}) how the convection matrix $\mathbb{N}$ can also be put in the Kronecker product form. It will be seen that this structure is useful for efficient matrix-vector product computations.

\subsection{Picard's method}
We use Picard's method to solve the nonlinear equation \cref{eq:ns_sys_all}. Picard's method is a fixed-point iteration. Let $\bm{u}^{(i)}$, $\bm{p}^{(i)}$ be the approximate solutions at the $i$th step. Each Picard step entails solving a large linear system
\begin{equation}
\label{eq:ns_sys_all_iter}
    \begin{pmatrix}
    \mathbb{F}(\bm{u}^{(i-1)})+\mathbb{C} & \mathbb{B}^T\\
    \mathbb{B} & 0
    \end{pmatrix}
    \begin{pmatrix} \bm{u}^{(i)}\\ \bm{p}^{(i)} \end{pmatrix} =
    \begin{pmatrix} \bm{f}^u\\ \bm{f}^p \end{pmatrix}.
\end{equation}
Instead of \cref{eq:ns_sys_all_iter}, one can equivalently solve the corresponding residual equation for a correction of the solution. Let $\bm{u}^{(i)}=\bm{u}^{(i-1)}+\delta\bm{u}^{(i)}$, $\bm{p}^{(i)}=\bm{p}^{(i-1)}+\delta\bm{p}^{(i)}$. Then $\delta\bm{u}^{(i)}$ and $\delta\bm{p}^{(i)}$ satisfy
\begin{equation}
\label{eq:ns_sys_all_res}
    \begin{pmatrix}
    \mathbb{F}(\bm{u}^{(i-1)})+\mathbb{C} & \mathbb{B}^T\\
    \mathbb{B} & 0
    \end{pmatrix}
    \begin{pmatrix} \delta\bm{u}^{(i)} \\ \delta\bm{p}^{(i)} \end{pmatrix} =
    \begin{pmatrix} \bm{r}^{u,(i-1)}\\ \bm{r}^{p,(i-1)} \end{pmatrix},
\end{equation}
where the nonlinear residual is
\begin{equation}
    \bm{r}^{(i)} =
    \begin{pmatrix} \bm{r}^{u,(i)}\\ \bm{r}^{p,(i)} \end{pmatrix} = 
    \begin{pmatrix} \bm{f}^u\\ \bm{f}^p \end{pmatrix} - 
    \begin{pmatrix}
    \mathbb{F}(\bm{u}^{(i)})+\mathbb{C} & \mathbb{B}^T\\
    \mathbb{B} & 0
    \end{pmatrix}
    \begin{pmatrix} \bm{u}^{(i)} \\ \bm{p}^{(i)} \end{pmatrix}.
\end{equation}
Let $\bm{f}$ denote the right-hand side of \cref{eq:ns_sys_all}. The complete algorithm is summarized in \cref{alg:picard}. The initial iterates $\bm{u}^{(0)}$, $\bm{p}^{(0)}$ are obtained as the solution of a Stokes problem, for which in \cref{eq:ns_sys_all_iter} the convection matrix $\mathbb{N}$ is set to zero.

\begin{algorithm}
\caption{Picard's method}
\label{alg:picard}
\begin{algorithmic}[1]
    \STATE{Solve Stokes problem for initial $\bm{u}^{(0)}$, $\bm{p}^{(0)}$, update convection matrix $\mathbb{N}(\bm{u}^{(0)})$, and compute nonlinear residual $\bm{r}^{(0)}$. $i=0$.}
    \WHILE{$\|\bm{r}^{(i)}\|_2>tol_\text{picard}*\|\bm{f}\|_2$ and $i<maxit$}
	\STATE{$i=i+1$}
	\STATE{Solve linear system \cref{eq:ns_sys_all_res} for $\delta\bm{u}^{(i)}$, $\delta\bm{p}^{(i)}$} 
	\STATE{Update solution $\bm{u}^{(i)}$, $\bm{p}^{(i)}$}
	\STATE{Update convection matrix $\mathbb{N}(\bm{u}^{(i)})$}
	\STATE{Compute nonlinear residual $\bm{r}^{(i)}$}
    \ENDWHILE
    \RETURN{$\bm{u}^{(i)}$, $\bm{p}^{(i)}$}
    \end{algorithmic}
\end{algorithm}

%-------------------------------------------------------------
\section{Low-rank approximation}
\label{sec:low_rank}

In this section we discuss low-rank approximation techniques and how they can be used with iterative solvers. The computational cost of solving \cref{eq:ns_sys_all_res} at each Picard step is high due to the large problem size $n_t n_\xi (n_u+n_p)$, especially when large numbers of spatial grid points or time steps are used to achieve high-resolution solution. We will address this using low-rank tensor approximations to the solution vectors $\bm{u}$ and $\bm{p}$. We will develop efficient iterative solvers and preconditioners where the solution is approximated using a compressed data representation in order to greatly reduce memory requirements and computational effort. The idea is to represent the iterates in an approximate Krylov subspace method in a low-rank tensor format. The basic operations associated with the low-rank format are much cheaper, and as the Krylov subspace method converges it constructs a sequence of low-rank approximations to the solution of the system.

\subsection{Tensor train decomposition}
A tensor $\underline{z}\in\mathbb{R}^{n_1\times\cdots\times n_d}$ is a multidimensional array with entries $\underline{z}(i_1,\ldots,i_d)$, where $i_l=1,\ldots,n_l$, $l=1,\ldots,d$. The solution coefficients in \cref{eq:sg_sol} can be represented in the form of three-dimensional $n_t\times n_\xi\times n_x$ tensors $\underline{u}$ (where $n_x=n_u$) and $\underline{p}$ ($n_x=n_p$), such that $\underline{u}(k,s,j)=u^k_{js}$ and $\underline{p}(k,s,j)=p^k_{js}$. Equivalently, such tensors can be represented in vector format, where the vector version $\bm{u}$ and $\bm{p}$ are specified using the vectorization operation
\begin{equation}
    \bm{u}=\text{vec}(\underline{u}) \,\,\Leftrightarrow\,\, \bm{u}(\overline{i_1i_2i_3})=\underline{u}(i_1,i_2,i_3)
\end{equation}
where $\overline{i_1i_2i_3}=i_3+(i_2-1)n_x+(i_1-1)n_\xi n_x$, and $\bm{p}=\text{vec}(\underline{p})$ in a similar manner. In an iterative solver for the system \cref{eq:ns_sys_all_res}, any iterate $\bm{z}$ can be equivalently represented as a three-dimensional tensor $\underline{z}\in\mathbb{R}^{n_t\times n_\xi\times n_x}$. In the sequel we use vector $\bm{z}$ and tensor $\underline{z}$ interchangebly. The tensor train decomposition \cite{Os11} is a compressed low-rank representation to approximate a given tensor and efficiently perform tensor operations. Specifically, the tensor train format of $\underline{z}$ is defined as
\begin{equation}
\label{eq:z_tt}
    \underline{z}(i_1,i_2,i_3) \approx \sum_{\alpha_1,\alpha_2} \underline{z}^{(1)}(i_1,\alpha_1) \underline{z}^{(2)}(\alpha_1,i_2,\alpha_2) \underline{z}^{(3)}(\alpha_2,i_3),
\end{equation}
where $\underline{z}^{(1)}\in\mathbb{R}^{n_t\times\kappa_1}$, $\underline{z}^{(2)}\in\mathbb{R}^{\kappa_1\times n_\xi\times\kappa_2}$, $\underline{z}^{(3)}\in\mathbb{R}^{\kappa_2\times n_x}$ are the tensor train cores, and $\kappa_1$ and $\kappa_2$ are called the tensor train ranks. It is easy to see that if $\kappa_1,\kappa_2\approx\kappa$ and $\kappa$ is small, the memory cost to store $\underline{z}$ is reduced from $O(n_tn_\xi n_x)$ to $O((n_t+n_\xi\kappa+n_x)\kappa)$. 

The tensor train decomposition allows efficient basic operations on tensors. Most importantly, matrix-vector products can be computed much less expensively if the vector $\bm{z}$ is in the tensor train format. For $\underline{z}$ as in \cref{eq:z_tt}, the vector $\bm{z}$ has an equivalent Kronecker product form \cite{DoSa14}
\begin{equation}
\label{eq:z_kron}
    \bm{z} = \text{vec}(\underline{z}) = \sum_{\alpha_1,\alpha_2} z^{(1)}_{\alpha_1}\otimes z^{(2)}_{\alpha_1,\alpha_2}\otimes z^{(3)}_{\alpha_2},
\end{equation}
where in the right-hand side $z^{(1)}_{\alpha_1}$, $z^{(2)}_{\alpha_1,\alpha_2}$, and $z^{(3)}_{\alpha_2}$ are vectors of length $n_t$, $n_\xi$, and $n_x$, respectively, obtained by fixing the indices $\alpha_1$ and $\alpha_2$ in $\underline{z}^{(1)}$, $\underline{z}^{(2)}$, and $\underline{z}^{(3)}$. Then for any matrix $\mathbb{X}= X^{(1)}\otimes X^{(2)}\otimes X^{(3)}$, such as the blocks in \cref{eq:ns_sys_all_res},
\begin{equation}
\label{eq:Xz}
    \mathbb{X}\bm{z} = \sum_{\alpha_1,\alpha_2} (X^{(1)}z^{(1)}_{\alpha_1}) \otimes (X^{(2)}z^{(2)}_{\alpha_1,\alpha_2}) \otimes (X^{(3)}z^{(3)}_{\alpha_2}).
\end{equation}
The product is also in tensor train format with the same ranks as in $\bm{z}$ (of the right-hand side of \cref{eq:z_tt}), and it only requires matrix-vector products for each component of $\mathbb{X}$. From left to right in the Kronecker products, the component matrices from \cref{eq:F} are sparse with numbers of nonzeros proportional to $n_t$, $n_\xi$, and $n_x$, respectively, and the computational cost is thus reduced from $O(n_tn_\xi n_x)$ to $O((n_t+n_\xi\kappa+n_x)\kappa)$.

Other vector computations, including additions and inner products, are also inexpensive with the tensor train format. One thing to note is that the additions of two vectors in tensor train format will tend to increase the ranks. This can be easily seen from \cref{eq:z_tt}, since the addition of two low-rank tensors end up with more terms for the summation on the right-hand side. An important operation for the tensor train format is a truncation (or rounding) operation, used to reduce the ranks for tensors that are already in the tensor train format but have suboptimal high ranks. For a given tensor $\underline{z}$ as in \cref{eq:z_tt}, the truncation operation $\mathcal{T}$ with tolerance $\epsilon$ computes
\begin{equation}
    \tilde{\underline{z}} = \mathcal{T}_\epsilon(\underline{z}),
\end{equation}
such that $\tilde{\underline{z}}$ has smaller ranks than $\underline{z}$ and satisfies the relative error
\begin{equation}
    \|\tilde{\underline{z}}-\underline{z}\|_F/\|\underline{z}\|_F\leq\epsilon.
\end{equation}
(Note that $\|\underline{z}\|_F=\|\bm{z}\|_2$.) The truncation operator is based on the TT-SVD algorithm \cite{Os11}, given in \cref{alg:tt_svd}, which is used to compute a low-rank tensor train approximation for a full tensor $\underline{z}\in\mathbb{R}^{n_1\times\cdots\times n_d}$. In the algorithm, a sequence of singular value decompositions (SVDs) is computed for the so-called unfolding matrix $Z$, obtained by reshaping the entries of a tensor into a two-dimensional array. Terms corresponding to small singular values are dropped such that an error $E$ in the truncated SVD satisfies $\|E\|_F\leq\delta_j$, $j=1,\ldots,d-1$ (see line 4 of \cref{alg:tt_svd}). It was shown in \cite{Os11} that the algorithm produces a tensor train $\tilde{\underline{z}}$ that satisfies
\begin{equation}
    \| \underline{z}-\tilde{\underline{z}} \|_F \leq \Big(\sum_{k=1}^{d-1}\delta_k^2 \Big)^{1/2}.
\end{equation}
Thus, one can choose $\delta_1=\cdots=\delta_{d-1}=\epsilon \|\underline{z}\|_F /\sqrt{d-1}$ to make the relative error $\|\tilde{\underline{z}}-\underline{z}\|_F/\|\underline{z}\|_F\leq\epsilon$. Note the algorithm is costly since it requires SVDs on matrices $Z\in\mathbb{R}^{\kappa_{j-1}n_j\times n_{j+1}\cdots n_d}$. However, when the tensor $\underline{z}$ is already in the tensor train format, the computation can be greatly simplified, and only SVDs on the much smaller tensor train cores are needed. In this case, the cost of the truncation operation is $O(dn\kappa^3)$ if $n_1,\ldots,n_d\approx n$ and $\kappa_1,\ldots,\kappa_{d-1}\approx \kappa$. We refer to \cite{Os11} for more details. In the numerical experiments, we use TT-Toolbox \cite{tt} for tensor train computations.
\begin{algorithm}
\caption{TT-SVD}
\label{alg:tt_svd}
\begin{algorithmic}[1]
    \STATE{Let $Z=\underline{z}$. Set truncation parameters $\{\delta_j\}$. $\kappa_0=1$.}
    \FOR{$j=1,\ldots,d-1$}
        \STATE{$Z\leftarrow \text{reshape}(Z,[\kappa_{j-1}n_j,n_{j+1}\cdots n_d])$}
        \STATE{Compute truncated SVD $Z=U\Sigma V^T+E$, $\|E\|_F\leq\delta_j$, $\kappa_j=\text{rank}(\Sigma)$}
        \STATE{New core $\tilde{\underline{z}}^{(j)}\leftarrow\text{reshape}(U,[\kappa_{j-1},n_j,\kappa_j])$}
        \STATE{Update $Z\leftarrow\Sigma V^T$}
    \ENDFOR
    \STATE{New core $\tilde{\underline{z}}^{(d)}\leftarrow Z$}
    \RETURN $\tilde{\underline{z}}$ in tensor train format with cores $\{\tilde{\underline{z}}^{(j)}\}$
    \end{algorithmic}
\end{algorithm}

\subsection{Low-rank solver}
The tensor train decomposition offers efficient tensor operations and we use it in iterative solvers to reduce the computational costs. The all-at-once system \cref{eq:ns_sys_all_res} to be solved at each step of Picard's method is nonsymmetric. We use a right-preconditioned GMRES method to solve the system. The complete algorithm for solving $\mathscr{L}\bm{z}=\bm{b}$ is summarized in \cref{alg:gmres}. The preconditioner $\mathscr{P}^{-1}$ entails an inner iterative process and is not fixed for each GMRES iteration, and therefore a variant of the flexible GMRES method (see, e.g., \cite{Sa03}) is used. As discussed above, all the iterates in the algorithm are represented in the tensor train format for efficient computations, and a truncation operation with tolerance $\epsilon_\text{gmres}$ is used to compress the tensor train ranks so that they stay small relative to the problem size. It should be noted that since the quantities are truncated, the Arnoldi vectors $\{\bm{v}_i\}$ do not form orthogonal basis for the Krylov subspace, and thus this is not a true GMRES computation. When the algorithm is used for solving  \cref{eq:ns_sys_all_res}, the truncation operator is applied to quantities associated with the two tensor trains $\delta\bm{u}^{(i)}$ and $\delta\bm{p}^{(i)}$ separately. In \cref{sec:prec}, we construct effective preconditioners for the system \cref{eq:ns_sys_all_res}.

\begin{algorithm}
\caption{Low-rank GMRES method}
\label{alg:gmres}
\begin{algorithmic}[1]
    \STATE{Choose initial $\bm{z}_0$, compute $\bm{s}_0=\mathcal{T}_{\epsilon_\text{gmres}}(\bm{b}-\mathscr{L}\bm{z}_0)$, $\beta=\|\bm{s}_0\|_2$, and $\bm{v}_1=\bm{s}_0/\beta$. Let $k=0$.}
    \WHILE{$\|\bm{s}_k\|_2>tol_\text{gmres}*\|\bm{b}\|_2$ and $k<maxit$}
        \STATE{$k = k+1$}
	\STATE{Compute $\hat{\bm{v}}_k=\mathscr{P}^{-1}\bm{v}_k$}
	\STATE{Compute $\bm{x}=\mathcal{T}_{\epsilon_\text{gmres}}(\mathscr{L}\hat{\bm{v}}_k$)}
	\FOR{$i=1,\ldots,k$}
	    \STATE{$h_{ik}=\bm{x}^T\bm{v}_i$}
	    \STATE{$\bm{x}=\bm{x}-h_{ik}\bm{v}_i$}
	\ENDFOR
	\STATE{$h_{k+1,k}=\|\bm{x}\|_2$, $\bm{v}_{k+1}=\mathcal{T}_{\epsilon_\text{gmres}}(\bm{x}/h_{k+1,k})$}
	\STATE{Define $\hat{V}_k=[\hat{\bm{v}}_1,\ldots,\hat{\bm{v}}_k]$ and $\bar{H}\in\mathbb{R}^{(k+1)\times k}$ with $\bar{H}_{ij}=h_{ij}$}
	\STATE{Compute $y_k=\text{argmin}_y\|\beta e_1-\bar{H}y\|_2$, where $e_1=[1,0,\ldots,0]^T$}
	\STATE{Compute $\bm{z}_k=\mathcal{T}_{\epsilon_\text{gmres}}(\bm{z}_0+\hat{V}_ky_k)$}
	\STATE{Compute $\bm{s}_k=\mathcal{T}_{\epsilon_\text{gmres}}(\bm{b}-\mathscr{L}\bm{z}_k)$}
    \ENDWHILE
    \RETURN $\bm{z}_k$
    \end{algorithmic}
\end{algorithm}

We also use the tensor train decomposition to construct a more efficient variant of \cref{alg:picard}. In particular, the updated solutions $\bm{u}^{(i)}$ and $\bm{p}^{(i)}$ in line 5 are truncated, with a tolerance ${\epsilon_\text{soln}}$, so that
\begin{equation}
\label{eq:eps_soln}
    \bm{u}^{(i)}=\mathcal{T}_{\epsilon_\text{soln}}(\bm{u}^{(i-1)}+\delta\bm{u}^{(i)}),\quad\bm{p}^{(i)}=\mathcal{T}_{\epsilon_\text{soln}}(\bm{p}^{(i-1)}+\delta\bm{p}^{(i)}).
\end{equation}
Another truncation operation with $\epsilon_\text{gmres}$ is applied to compress the ranks of the nonlinear residual $\bm{r}^{(i)}$ in line 7. We will use this truncated version of \cref{alg:picard} in numerical experiments; choices of the truncation tolerances will be specified in \cref{sec:num_test}.

\subsection{Convection matrix}
\label{sec:stoch_conv}
We now show that in \cref{eq:F} if the velocity $\bm{u}$ is in the tensor train format, the convection matrix $\mathbb{N}(\bm{u})$ can be represented as a sum of Kronecker products of matrices \cite{BeDo17}, which allows efficient matrix-vector product computations as in \cref{eq:Xz}. Assume the coefficient tensor in \cref{eq:sg_sol} is approximated by a tensor train decomposition,
\begin{equation}
    u_{jl}^k=\underline{u}(k,l,j) = \sum_{\alpha_1,\alpha_2} \underline{u}^{(1)}(k,\alpha_1) \underline{u}^{(2)}(\alpha_1,l,\alpha_2) \underline{u}^{(3)}(\alpha_2,j).
\end{equation}
Note that the entries of $\bm{N}(\vec{u}_{h,l}^k)$ are linear in $\vec{u}_{h,l}^k$ and
\begin{equation}
    \vec{u}_{h,l}^k = \sum_j u_{jl}^k\vec{\phi}_j(x) = \sum_{\alpha_1,\alpha_2} \underline{u}^{(1)}(k,\alpha_1) \underline{u}^{(2)}(\alpha_1,l,\alpha_2) (\sum_j \underline{u}^{(3)}(\alpha_2,j)\vec{\phi}_j(x)).
\end{equation}
Let $\vec{u}_{\alpha_2}^{(3)}=\sum_j \underline{u}^{(3)}(\alpha_2,j)\vec{\phi}_j(x)$. Then the $k$th diagonal block of $\mathbb{N}(\bm{u})$ is
\begin{equation}
    \sum_{l=1}^{n_\xi} (H_l\otimes\bm{N}(\vec{u}_{h,l}^k)) = \sum_{\alpha_1,\alpha_2} \underline{u}^{(1)}(k,\alpha_1) \sum_{l=1}^{n_\xi}(\underline{u}^{(2)}(\alpha_1,l,\alpha_2) H_l) \otimes \bm{N}(\vec{u}_{\alpha_2}^{(3)}).
\end{equation}
The convection matrix $\mathbb{N}(\bm{u})$ can be expressed as
\begin{equation}
\label{eq:N_kron}
    \mathbb{N}(\bm{u}) = \sum_{\alpha_1,\alpha_2} \text{diag}(u_{\alpha_1}^{(1)}) \otimes \sum_{l=1}^{n_\xi}(\underline{u}^{(2)}(\alpha_1,l,\alpha_2) H_l) \otimes \bm{N}(\vec{u}_{\alpha_2}^{(3)}).
\end{equation}
Here $u^{(1)}_{\alpha_1}$ is a vector obtained by fixing the index $\alpha_1$ in $\underline{u}^{(1)}$, and $\text{diag}(u^{(1)}_{\alpha_1})$ is a diagonal matrix with $u^{(1)}_{\alpha_1}$ on the diagonal. The result is a sum of Kronecker products of three smaller matrices. Such a representation can be constructed for any iterate $\bm{u}^{(i)}$ in the tensor train format. 

Given the number of terms in the summation in the right-hand side of \cref{eq:N_kron}, the matrix-vector product with $\mathbb{N}$ will result in a dramatic tensor train rank increase, from $\kappa$ to $\kappa^2$. Unless $\kappa$ is very small, a tensor train with rank $\kappa^2$ will require too much memory and also be expensive to work with. To overcome this difficulty, when solving the all-at-once system \cref{eq:ns_sys_all_res}, we use a low-rank approximation of $\bm{u}^{(i)}$ to construct $\mathbb{N}(\bm{u}^{(i)})$. Specifically, let
\begin{equation}
\label{eq:eps_conv}
    \tilde{\bm{u}}^{(i)} = \mathcal{T}_{\epsilon_\text{conv}} (\bm{u}^{(i)})
\end{equation}
with some truncation tolerance $\epsilon_\text{conv}$. Since $\tilde{\bm{u}}^{(i)}$ has smaller ranks than $\bm{u}^{(i)}$, the approximate convection matrix $\mathbb{N}(\tilde{\bm{u}}^{(i)})$ contains a smaller number of terms in \cref{eq:N_kron}, and thus the rank increase becomes less significant when computing matrix-vector products with it. In other words, the linear system solved at each Picard step becomes
\begin{equation}
    \begin{pmatrix}
    \mathbb{F}(\tilde{\bm{u}}^{(i-1)})+\mathbb{C} & \mathbb{B}^T\\
    \mathbb{B} & 0
    \end{pmatrix}
    \begin{pmatrix} \delta\bm{u}^{(i)} \\ \delta\bm{p}^{(i)} \end{pmatrix} =
    \begin{pmatrix} \bm{r}^{u,(i-1)}\\ \bm{r}^{p,(i-1)} \end{pmatrix}.
\end{equation}
Note that the original $\bm{u}^{(i)}$ is still used for computing the nonlinear residual $\bm{r}^{(i)}$ in Picard's method.

%-------------------------------------------------------------
\section{Preconditioning} 
\label{sec:prec}

In this section we discuss preconditioning techniques for the all-at-once system \cref{eq:ns_sys_all_res} so that the Krylov subspace methods converge in a small number of iterations. To simplify the notation, we use $\bm{w}$ instead of $\bm{u}^{(i-1)}$, and the associated approximate solution at the $k$th time step is
\begin{equation}
    \vec{w}_{h}^k(x,\xi) = \sum_{l=1}^{n_\xi}\sum_{j=1}^{n_u} w_{jl}^k \vec{\phi}_j(x) \psi_l(\xi) = \sum_{l=1}^{n_\xi} \vec{w}_{h,l}^k(x)\psi_l(\xi)
\end{equation}
with $\vec{w}_{h,l}^k(x)=\sum_j w_{jl}^k\vec{\phi}_j(x)$. In the following the dependence on $\bm{w}$ in $\mathbb{F}(\bm{w})$ is omitted in most cases. We derive a preconditioner by extending ideas for more standard problems \cite{ElSi14}, starting with an ``idealized'' block triangular preconditioner
\begin{equation}
\label{eq:prec}
    \mathscr{P} = \begin{pmatrix}
    \mathbb{F}+\mathbb{C} & \mathbb{B}^T\\
    0 & -\mathbb{S}
    \end{pmatrix}.
\end{equation}
With this choice of preconditioner, the Schur complement is $\mathbb{S}=\mathbb{B}(\mathbb{F}+\mathbb{C})^{-1}\mathbb{B}^T$, and the idealized preconditioned system derived from a block factorization
\begin{equation}
    \begin{pmatrix} \mathbb{F}+\mathbb{C} & \mathbb{B}^T\\
    \mathbb{B} & 0 \end{pmatrix} \mathscr{P}^{-1} =
    \begin{pmatrix} \mathbb{I} & 0 \\
    \mathbb{B}\mathbb{F}^{-1} & \mathbb{I} \end{pmatrix}
\end{equation}
has eigenvalues equal to 1 and Jordan blocks of order 2. (Here $\mathbb{I}$ is an identity block.) Thus a right-preconditioned true GMRES method will converge in two iterations. However, the application of $\mathscr{P}^{-1}$ involves solving linear systems associated with $\mathbb{S}$ and $\mathbb{F}+\mathbb{C}$. These are too expensive for practical computation and to develop preconditioners we will construct inexpensive approximations to the linear solves. Specifically, we derive mean-based preconditioners that use results from the mean deterministic problem. Such preconditioners for the stochastic steady-state Navier--Stokes equations have been studied in \cite{PoSi12}. We generalize the techniques for the all-at-once formulation of the unsteady equations.

\subsection{Deterministic operator}
We review the techniques used for approximating the Schur complement in the deterministic case \cite{ElSi14}. The approximations are based on the fact that a commutator of the convection-diffusion operator with the divergence operator
\begin{equation}
\label{eq:commu_cont}
    \mathcal{E} = \nabla\cdot(-\nu\nabla^2+\vec{w}^k_{h,1}\cdot\nabla) - (-\nu\nabla^2+\vec{w}^k_{h,1}\cdot\nabla)_p\nabla\cdot
\end{equation}
is small under certain assumptions about smoothness and boundary conditions. The subscript $p$ means the operators are defined on the pressure space. For a discrete convection-diffusion operator $\bm{F}=\bm{A}_0+\bm{N}(\vec{w}_{h,1}^k)$ (which is part of the mean problem we discuss later), as defined in \cref{eq:mat_entries}, an approximation to the Schur complement $S = B\bm{F}^{-1}B^T$ is identified from a discrete analogue of \cref{eq:commu_cont},
\begin{equation}
\label{eq:commutator}
    E = (M_p^{-1}B)(\bm{M}^{-1}\bm{F}) - (M_p^{-1}F_p)(M_p^{-1}B) \approx 0,
\end{equation}
where the subscript $p$ means the corresponding matrices constructed on the discrete pressure space. \Cref{eq:commutator} leads to an approximation to the Schur complement matrix,
\begin{equation}
\label{eq:schur_approx}
    S = B\bm{F}^{-1}B^T \approx M_pF_p^{-1}B\bm{M}^{-1}B^T.
\end{equation}
The pressure convection-diffusion (PCD) preconditioner is constructed by replacing the mass matrices with approximations containing only their diagonal entries (denoted by a subscript $*$) in \cref{eq:schur_approx},
\begin{equation}
    S_\text{PCD}^{-1} = (B\bm{M_*}^{-1}B^T)^{-1}F_p M_{p*}^{-1}.
\end{equation}
The least-squares commutator (LSC) preconditioner avoids the construction of matrices on the pressure space, with the approximation to $F_p$,
\begin{equation}
\label{eq:lsc_approx}
    F_p \approx (B\bm{M}^{-1}\bm{F}\bm{M}^{-1}B^T)(B\bm{M}^{-1}B^T)^{-1}M_p
\end{equation}
(see \cite[section~9.2]{ElSi14} for a derivation). The LSC preconditioner is obtained by substituting $F_p$ in \cref{eq:schur_approx} and replacing the mass matrices with their diagonals,
\begin{equation}
    S_\text{LSC}^{-1} = (B\bm{M}_*^{-1}B^T)^{-1}(B\bm{M}_*^{-1}\bm{F}\bm{M}_*^{-1}B^T)(B\bm{M}_*^{-1}B^T)^{-1}.
\end{equation}
For both preconditioners, the only use of the matrices $\bm{F}$ and $F_p$ is through matrix-vector products with them.

\subsection{Approximations to $\mathbb{S}^{-1}$}
The Schur complement $\mathbb{S}$ involves $(\mathbb{F}+\mathbb{C})^{-1}$ and is impractical to work with. For our stochastic unsteady problem, we consider mean-based preconditioners that use approximations to the Schur complement matrix
\begin{equation}
    \mathbb{S}_0 = \mathbb{B}(\mathbb{F}_0+\mathbb{C})^{-1}\mathbb{B}^T,
\end{equation}
where the ``mean'' matrix $\mathbb{F}_0$ is block-diagonal with $\mathbb{F}_0^k$ as the $k$th diagonal block, and
\begin{equation}
    \mathbb{F}_0^k = \tau^{-1}(I_{n_\xi}\otimes\bm{M}) + I_{n_\xi}\otimes\bm{A}_0 + I_{n_\xi}\otimes\bm{N}(\vec{w}_{h,1}^k).
\end{equation}
This corresponds to taking only the first term in the two summations on the right-hand side of \cref{eq:F^k}. Since the gPC basis functions are orthonormal with $\langle\psi_r\psi_s\rangle=\delta_{rs}$ and $\psi_1\equiv 1$, it follows $\langle\psi_s\rangle=\delta_{1s}$, and $G_0=H_1=I_{n_\xi}$. 
The matrices $\bm{A}_0$ and $\bm{N}(\vec{w}_{h,1}^k)$ are constructed from the mean of $\nu$ and $\vec{w}_h^k$,
\begin{equation}
    \langle\nu\rangle=\nu_0,\,\,\langle\vec{w}_h^k\rangle=\sum\nolimits_jw_{j1}^k\vec{\phi}_j(x)=\vec{w}^k_{h,1}.
\end{equation}
The matrix $\mathbb{F}_0^k$ can be expressed as $I_{n_\xi}\otimes(\tau^{-1}\bm{M}+\bm{A}_0+\bm{N}(\vec{w}_{h,1}^k))$ and this enables use of approximations associated with a deterministic problem. Now, similarly define $\mathbb{F}_{p,0}$ on the pressure space, with
\begin{equation}
    \mathbb{F}_{p,0}^k = \tau^{-1}(I_{n_\xi}\otimes M_p) + I_{n_\xi}\otimes A_{p,0} + I_{n_\xi}\otimes N_p(\vec{w}_{h,1}^k).
\end{equation}
Let $\mathbb{M}=I_{n_t}\otimes I_{n_\xi}\otimes\bm{M}$ and $\mathbb{M}_p=I_{n_t}\otimes I_{n_\xi}\otimes M_p$. Assuming the validity of \cref{eq:commutator} it is easy to check that 
\begin{equation}
\label{eq:stoch_commu_1}
    \mathbb{M}_p^{-1}\mathbb{B}\mathbb{M}^{-1}\mathbb{F}_0 - \mathbb{M}_p^{-1}\mathbb{F}_{p,0}\mathbb{M}_p^{-1}\mathbb{B} \approx 0.
\end{equation}
On the other hand, let $\mathbb{C}_p=-\tau^{-1}C_{n_t}\otimes I_{n_\xi}\otimes M_p$, so that $\mathbb{C}$ satisfies
\begin{equation}
\label{eq:stoch_commu_2}
    \mathbb{M}_p^{-1}\mathbb{B}\mathbb{M}^{-1}\mathbb{C} - \mathbb{M}_p^{-1}\mathbb{C}_p\mathbb{M}_p^{-1}\mathbb{B} = 0.
\end{equation}
Combining \cref{eq:stoch_commu_1} and \cref{eq:stoch_commu_2} gives an approximation to $\mathbb{S}_0$,
\begin{equation}
    \mathbb{S}_0 = \mathbb{B}(\mathbb{F}_0+\mathbb{C})^{-1}\mathbb{B}^T \approx \mathbb{M}_p(\mathbb{F}_{p,0}+\mathbb{C}_p)^{-1}\mathbb{B}\mathbb{M}^{-1}\mathbb{B}^T.
\end{equation}
Then the mean-based PCD preconditioner is given as
\begin{equation}
\label{eq:mean_pcd}
    \mathbb{S}_{\text{PCD},0}^{-1} = (\mathbb{B}\mathbb{M}_*^{-1}\mathbb{B}^T)^{-1}(\mathbb{F}_{p,0}+\mathbb{C}_p)\mathbb{M}_{p*}^{-1},
\end{equation}
where $\mathbb{M}_*=I_{n_t}\otimes I_{n_\xi}\otimes\bm{M}_*$ and $\mathbb{M}_{p*}=I_{n_t}\otimes I_{n_\xi}\otimes {M}_{p*}$. Similarly from \cref{eq:lsc_approx}, it holds that
\begin{equation}
    \mathbb{F}_{p,0}+\mathbb{C}_p \approx (\mathbb{B}\mathbb{M}^{-1}(\mathbb{F}_0+\mathbb{C})\mathbb{M}^{-1}\mathbb{B}^T)(\mathbb{B}\mathbb{M}^{-1}\mathbb{B})^{-1}\mathbb{M}_p.
\end{equation}
Substituting $\mathbb{F}_{p,0}+\mathbb{C}_p$ in \cref{eq:mean_pcd} and replacement of the mass matrices with their diagonals gives the mean-based LSC preconditioner
\begin{equation}
\label{eq:mean_lsc}
    \mathbb{S}_{\text{LSC},0}^{-1} = (\mathbb{B}\mathbb{M}_*^{-1}\mathbb{B}^T)^{-1}(\mathbb{B}\mathbb{M}_*^{-1}(\mathbb{F}_0+\mathbb{C})\mathbb{M}_*^{-1}\mathbb{B}^T)(\mathbb{B}\mathbb{M}_*^{-1}\mathbb{B}^T)^{-1}.
\end{equation}
The two mean-based preconditioners in \cref{eq:mean_pcd,eq:mean_lsc} have the same form as for the deterministic problem, except that there is an extra term $\mathbb{C}$ or $\mathbb{C}_p$ from the all-at-once formulation. Computations associated with the two approximations to the Schur complement are also inexpensive. For example, $(\mathbb{B}\mathbb{M}_*^{-1}\mathbb{B}^T)^{-1}=I_{n_t}\otimes I_{n_\xi}\otimes (B\bm{M}_*^{-1}B^T)^{-1}$, and this only requires solving a system with $B\bm{M}_*^{-1}B^T$ a discrete Laplacian. Multiplications with the mean matrix $\mathbb{F}_0+\mathbb{C}$ are reduced to its components (see \cref{eq:Xz}),
\begin{equation}
    \mathbb{F}_0+\mathbb{C} = \tau^{-1}(I_{n_t}\otimes I_{n_\xi}\otimes\bm{M}) + I_{n_t}\otimes I_{n_\xi}\otimes\bm{A}_0 + \mathbb{N}_0 - \tau^{-1}(C_{n_t}\otimes I_{n_\xi}\otimes\bm{M}).
\end{equation}
The matrix $\mathbb{N}_0$ is block-diagonal with $\mathbb{N}_0^k=I_{n_\xi}\otimes\bm{N}(\vec{w}_{h,1}^k)$ and can be expressed as a sum of Kronecker products of matrices as discussed in \cref{sec:stoch_conv},
\begin{equation}
    \mathbb{N}_0(\bm{w}) = \sum_{\alpha_1,\alpha_2} \text{diag}(w_{\alpha_1}^{(1)}) \otimes (\underline{w}^{(2)}(\alpha_1,1,\alpha_2) I_{n_\xi}) \otimes \bm{N}(\vec{w}_{\alpha_2}^{(3)}).
\end{equation}

\subsection{System solve with $\mathbb{F}+\mathbb{C}$}
The application of the preconditioner $\mathscr{P}^{-1}$ in \cref{eq:prec} also involves solving a linear system associated with the (1,1) block $\mathbb{F}+\mathbb{C}$. For approximation, we replace it with the mean matrix $\mathbb{F}_0+\mathbb{C}$, and solve a system of the form
\begin{equation}
\label{eq:11_sys}
    (\mathbb{F}_0+\mathbb{C}) \bm{v} = \bm{y}.
\end{equation}
For such a system it is easy to compute matrix-vector products and we again use a low-rank GMRES method for solving the system. This inner GMRES solver is preconditioned with
\begin{equation}
\label{eq:11_prec}
    \mathscr{M} = (I_{n_t}-C_{n_t})\otimes I_{n_\xi}\otimes (\tau^{-1}\bm{M}+\bm{A}_0+\bm{N}(\vec{w}_{h,1}^{\text{avg}})),
\end{equation}
where $\vec{w}_{h,1}^{\text{avg}}$ is the average of $\vec{w}_{h,1}^{k}$ over all time steps. For small time step $\tau$, the contribution from the mass matrix, $\tau^{-1}\bm{M}$, becomes dominant and $\mathscr{M}$ forms a good approximation to the coefficient matrix $\mathbb{F}_0+\mathbb{C}$. The application of $\mathscr{M}^{-1}$ is also conveniently reduced to computations associated with smaller matrices. We note that \cref{eq:11_sys} need not be solved accurately. In particular, with a stopping criterion $\|\bm{y}-(\mathbb{F}_0+\mathbb{C}) \bm{v}\|_2 \leq tol\|\bm{y}\|_2$, a relatively large stopping tolerance, e.g., $tol=10^{-1}$, will suffice for the mean-based preconditioner $\mathscr{P}$ to be effective.
\begin{remark}
For systems like \cref{eq:11_sys}, a block diagonal preconditioner ($\mathscr{M}=\mathbb{F}_0$) was studied in \cite{McWa16}, where it was shown that preconditioned GMRES converges very slowly before a sharp drop in the residual occurs when the number of iterations reaches $n_t$, which is equal to the number of diagonal blocks. In numerical experiments, we found that the preconditioner in \cref{eq:11_prec} is more effective than a block diagonal one, for which performance deteriorates as $\tau$ becomes smaller.
\end{remark}

%-------------------------------------------------------------
\section{Numerical experiments}
\label{sec:num_test}

\subsection{Benchmark problem}
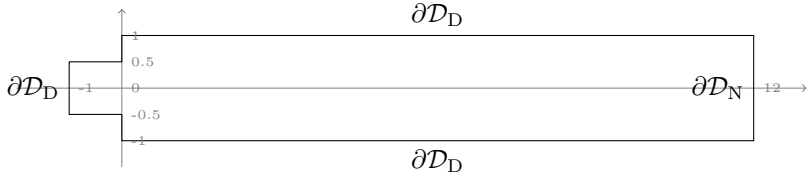
\begin{figure}[htbp]
\label{fig:domain}
\centering
\begin{tikzpicture}[scale=1.4]
    \draw [help lines, ->] (-1,0) -- (6.5,0);
    \draw [help lines, ->] (0,-0.75) -- (0,0.75);
    \draw (-0.5,0.25) -- (0,0.25) -- (0,0.5) -- (6,0.5) -- (6,-0.5) -- (0,-0.5) -- (0,-0.25) -- (-0.5,-0.25) -- (-0.5,0.25);
    \node [left] at (6,0) {$\partial\mathcal{D}_\text{N}$};
    \node [left] at (-0.5,0) {$\partial\mathcal{D}_\text{D}$};
    \node [above] at (3,0.5) {$\partial\mathcal{D}_\text{D}$};
    \node [below] at (3,-0.5) {$\partial\mathcal{D}_\text{D}$};
    \node [right, gray] at (0,0) {\tiny{0}};
    \node [right, gray] at (0,0.25) {\tiny{0.5}};
    \node [right, gray] at (0,0.5) {\tiny{1}};
    \node [right, gray] at (0,-0.25) {\tiny{-0.5}};
    \node [right, gray] at (0,-0.5) {\tiny{-1}};
    \node [right, gray] at (-0.5,0) {\tiny{-1}};
    \node [right, gray] at (6,0) {\tiny{12}};
\end{tikzpicture}
\caption{Symmetric step domain with boundary conditions.}
\end{figure}

Consider a flow around a symmetric step where the spatial domain $\mathcal{D}$ is a two-dimensional rectangular duct with a symmetric expansion (see \cref{fig:domain}). The Dirichlet inflow boundary conditions at $(-1,x_2)$, $|x_2|\leq0.5$ are deterministic and time-dependent, growing from zero to a steady parabolic profile,
\begin{equation}
    \vec{u}_\text{D}((-1,x_2),t) = \begin{pmatrix}1-4x_2^2 \\ 0 \end{pmatrix}(1-e^{-10t}).
\end{equation}
Neumann boundary conditions $\nu {\partial u_{x_1}}/{\partial x_1}=p$, ${\partial u_{x_2}}/{\partial x_1}=0$ are imposed at the outflow boundary $(12,x_2)$, $|x_2|\leq 1$, and no-flow conditions $\vec{u}=\vec{0}$ at the fixed walls $(x_1,\pm1)$, $0\leq x_1\leq 12$; $(x_1,\pm 0.5)$, $-1\leq x_1\leq 0$; $(0,x_2)$, $0.5\leq |x_2|\leq 1$. The initial conditions are zero everywhere for both $\vec{u}$ and $p$. The Taylor--Hood spatial discretization with biquadratic basis functions for the velocity space and bilinear basis functions for the pressure space is defined on a uniform grid of square elements with mesh size $h$, and it is constructed using the IFISS software package \cite{ifiss}.

The stochastic viscosity $\nu(x,\xi)$ is represented as a truncated KL expansion
\begin{equation}
    \nu(x,\xi) = \nu_0\Big(1.0 + \sigma\sum_{l=1}^{m} \sqrt{\beta_l}a_l(x)\xi_l\Big).
\end{equation}
The constants $\nu_0$ and $\nu_0\sigma$ represent the mean and the standard deviation of the stochastic field. We use an exponential covariance function $c(x,y)=\exp(-\|x-y\|_1/b)$, where $b$ is the correlation length. The pair $(\beta_l,a_l(x))$ is the $l$th largest eigenvalue and the corresponding eigenfunction of $c(x,y)$, satisfying
\begin{equation}
    \int_\mathcal{D} c(x,y)a_l(y)\text{d}y = \beta_la_l(x).
\end{equation}
This can be computed with a standard finite element method. The random variables $\{\xi_l\}_{l=1}^m$ are assumed to be independent and each of them uniformly distributed on the interval $[-\sqrt{3},\sqrt{3}]$, so they have zero means and unit variances. For the stochastic Galerkin method, the basis functions $\{\psi_r\}_{r=1}^{n_\xi}$ are $m$-dimensional Legendre polynomials, with total degrees bounded by $d_\psi$. Then the number of stochastic basis functions is $n_\xi=(m+d_\psi)!/(m!d_\psi !)$. In the numerical experiments, unless otherwise stated, the parameter values associated with the discrete problem are chosen as in \cref{table:parameters}. This gives a problem with dimensions $n_t=64$, $n_\xi=20$, $n_u=2992$, $n_p=461$, and $n_t n_\xi (n_u+n_p)=4419840$. All computations are done in MATLAB 9.4.0 (R2018a) on a desktop with 64 GB memory.
\begin{table}[hbtp]
\label{table:parameters}
\caption{Parameter values for numerical experiments.}
\centering
\begin{tabular}{c*{7}{|c}}
\hline
$\nu_0$ & $\sigma$ & $b$ & $m$ & $d_\psi$ & $t_f$ & $\tau$ & $h$\\ \hline
$1/50$ & 0.01 & 4.0 & 3 & 3 & 1.0 & $2^{-6}$ & $2^{-2}$\\
\hline
\end{tabular}
\end{table}

\subsection{Inexact Picard method}
\label{sec:inexact_picard}
The main computational cost associated with Picard's method is to solve an all-at-once system \cref{eq:ns_sys_all_res} at each step. In \cref{sec:low_rank} we discussed how to construct low-rank approximate solutions in tensor train format with much cheaper computations. To further reduce the cost, we adopt the idea of inexact Picard method \cite{Bi15}, where the linear systems are solved inexactly to save unnecessary computational work. Let \cref{eq:ns_sys_all_res} be denoted as $\mathscr{L}\bm{z}^{(i)}=\bm{r}^{(i-1)}$, and define the residual norm $\| \bm{s}_k \|_2 = \| \bm{r}^{(i-1)}-\mathscr{L}\bm{z}^{(i)}_k \|_2$ for an approximate solution $\bm{z}^{(i)}_k$. It was shown in \cite{Bi15} that if the stopping criterion for the linear solve (line 2 of \cref{alg:gmres}) is given as
\begin{equation}
\label{eq:inexact_solve}
    \| \bm{s}_k \|_2 \leq tol_\text{gmres} \| \bm{r}^{(i-1)} \|_2,
\end{equation}
then Picard's method converges as long as $tol_\text{gmres}<1$. This is especially helpful for our low-rank GMRES method. The best accuracy that the low-rank GMRES method can achieve is related to the truncation tolerance $\epsilon_\text{gmres}$ used in the algorithm (see \cref{fig:converge_gmres}). A relaxed stopping tolerance not only reduces the number of GMRES iterations, but it also allows use of larger truncation tolerances for tensor rank compressions, resulting in smaller ranks for the iterates and more efficient computations in the iterative solver. In the numerical tests, we set $tol_\text{gmres}=10^{-1}$ and $\epsilon_\text{gmres}=10^{-2}*tol_\text{gmres}=10^{-3}$. The same tolerances are used for solving the linear system \cref{eq:11_sys} required for the preconditioning operation. For the initial $\bm{u}^{(0)}$, $\bm{p}^{(0)}$, the Stokes problem is solved to satisfy $\| \bm{s}_k \|_2 \leq tol_\text{gmres} \| \bm{f} \|_2$ where $\bm{f}$ is the right-hand side of \cref{eq:ns_sys_all_iter}.

\subsection{Numerical results}
\begin{figure}[hbtp]
    \label{fig:converge}
    \begin{center}
    \subfloat[]{\label{fig:converge_gmres}\includegraphics[width=0.5\textwidth]{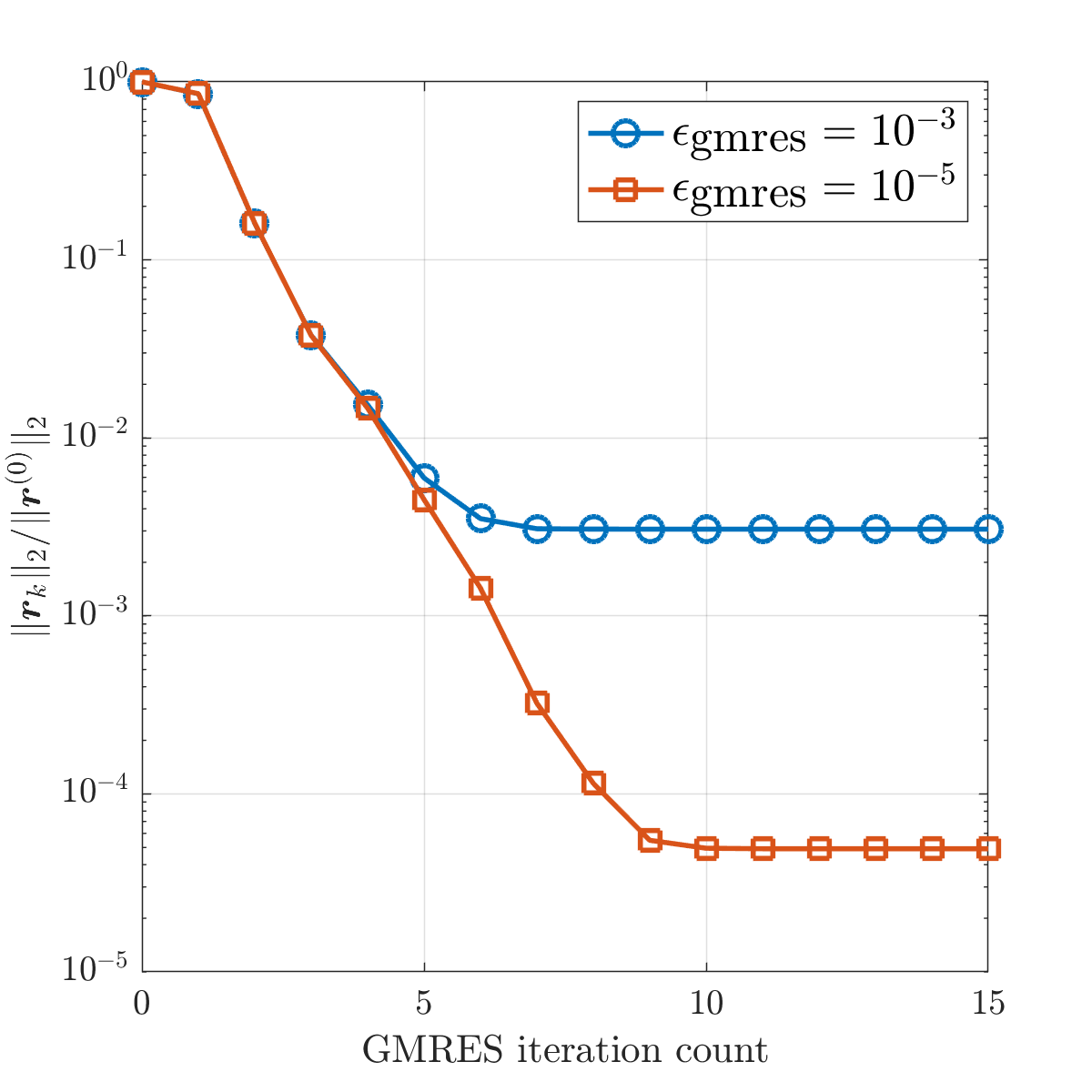}}
    \subfloat[]{\label{fig:converge_picard}\includegraphics[width=0.5\textwidth]{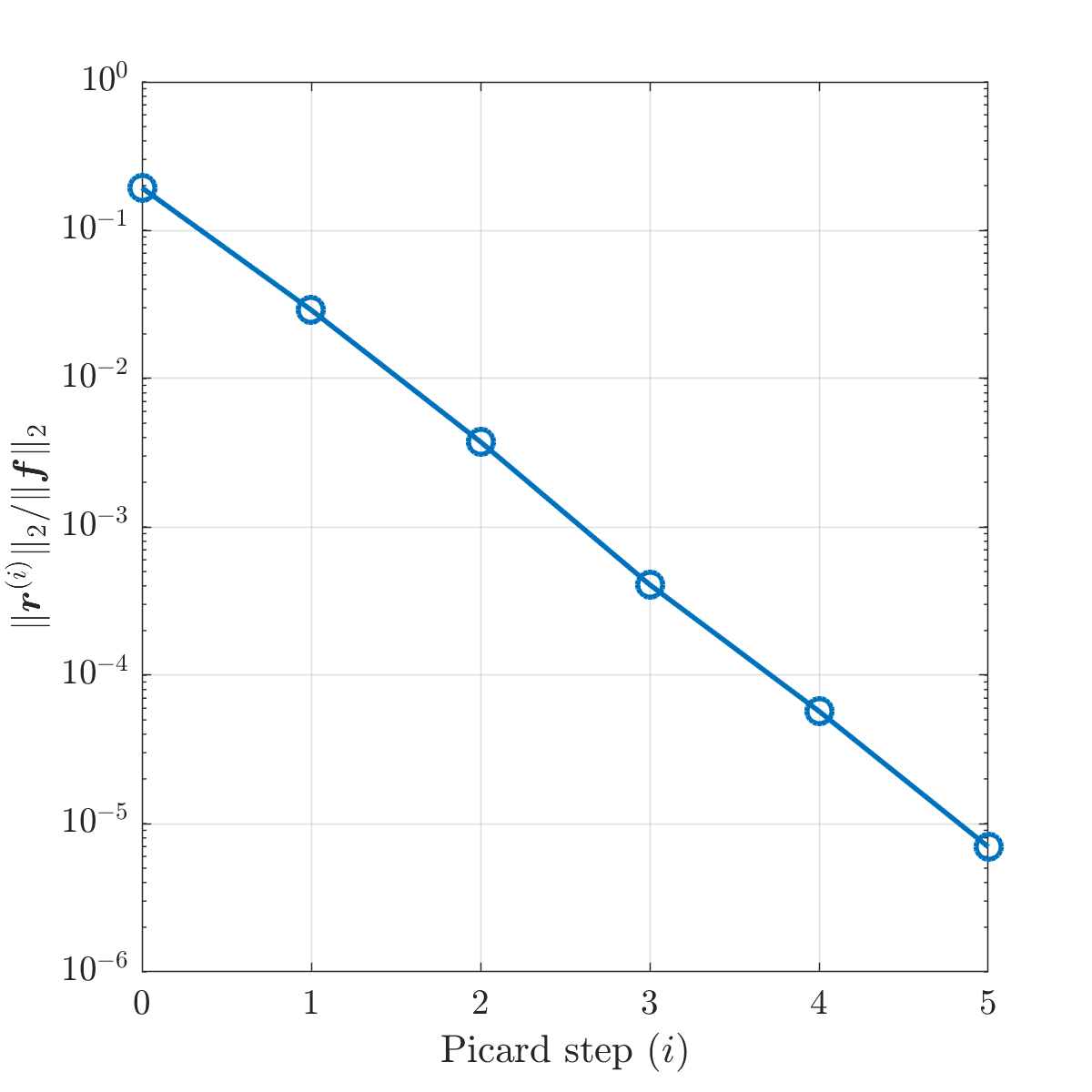}}
    \end{center}
    \caption{(a) Convergence of the low-rank GMRES method (at the first Picard step) with different truncation tolerances. (b) Convergence of the inexact Picard method with tolerances chosen as in \cref{table:tolerances}. LSC preconditioner is used.}
\end{figure}

\begin{figure}[hbtp]
    \label{fig:ranks_step}
    \begin{center}
    \subfloat[]{\label{fig:ranks_step_du}\includegraphics[width=0.5\textwidth]{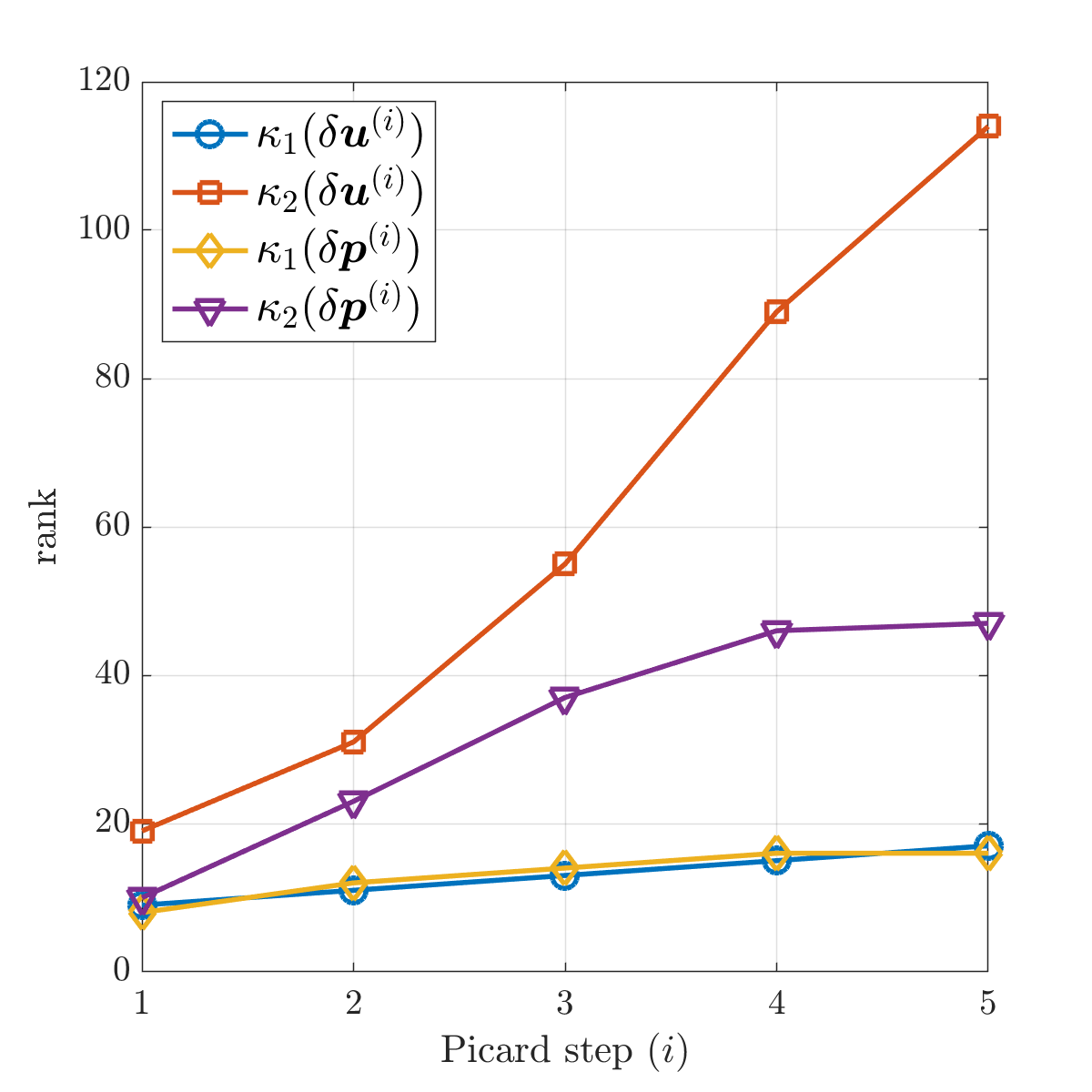}}
    \subfloat[]{\label{fig:ranks_step_u}\includegraphics[width=0.5\textwidth]{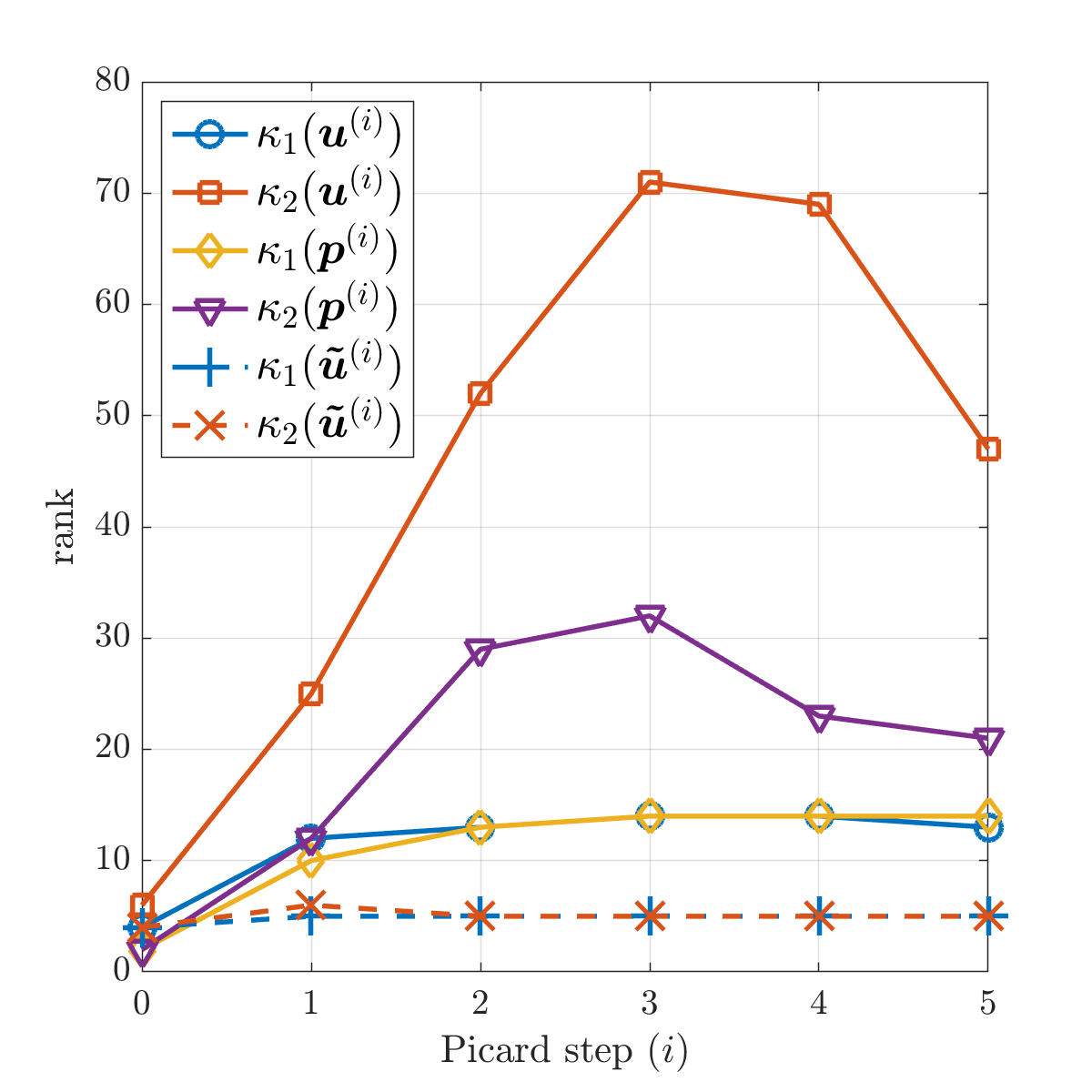}}
    \end{center}
    \caption{(a) Tensor train ranks of corrections $\delta\bm{u}^{(i)}$ and $\delta\bm{p}^{(i)}$. (b) Tensor train ranks of approximate solutions $\bm{u}^{(i)}$ and $\bm{p}^{(i)}$, and tensor train ranks of $\tilde{\bm{u}}^{(i)}$ for convection matrix. LSC preconditioner is used.}
\end{figure}

In the following, we examine the performance of the proposed low-rank algorithm in different settings. The choices of stopping and truncation tolerances are summarized in \cref{table:tolerances}. In \cref{alg:picard}, the stopping criterion for Picard's method is
\begin{equation}
    \|\bm{r}^{(i)}\|_2 \leq tol_\text{picard}*\|\bm{f}\|_2.
\end{equation}
We set $tol_\text{picard}=10^{-5}$. A small truncation tolerance $\epsilon_\text{soln}=10^{-7}$ is used to produce low-rank approximate solutions $\bm{u}^{(i)}$ and $\bm{p}^{(i)}$ in \cref{eq:eps_soln}. It is shown in \cref{fig:converge_picard} that, like the exact method, the inexact Picard method still exhibits a linear convergence rate. It takes 5 Picard steps to reach the required accuracy. \Cref{fig:ranks_step} shows the tensor train ranks $\kappa_1$ and $\kappa_2$ of the iterates at each Picard step. As the Picard iteration converges, the right-hand side of \cref{eq:inexact_solve} becomes smaller, and the corrections $\delta\bm{u}^{(i)}$ and $\delta\bm{p}^{(i)}$ computed from the low-rank GMRES method have increasing ranks. On the other hand, for the approximate solutions $\bm{u}^{(i)}$ and $\bm{p}^{(i)}$, their ranks drop to smaller values in the latter steps of the iteration. With a more stringent $tol_\text{picard}$, a smaller $\epsilon_\text{soln}$ is required and the approximate solutions have slightly higher ranks than those shown in \cref{fig:ranks_step_u}. Also shown in \cref{fig:ranks_step_u} are the tensor train ranks of $\tilde{\bm{u}}^{(i)}$ for constructing the approximate convection matrices using \cref{eq:eps_conv}. They have much smaller values than the ranks of $\bm{u}^{(i)}$.

\begin{table}[hbtp]
\label{table:tolerances}
\caption{Stopping and truncation tolerances.}
\centering
\begin{tabular}{c|r}
\hline
GMRES stopping tolerance & $tol_\text{gmres}=10^{-1}$ \\
Picard stopping tolerance & $tol_\text{picard}=10^{-5}$ \\
GMRES truncation tolerance & $\epsilon_\text{gmres}=10^{-3}$ \\
Truncation tolerance for solutions & $\epsilon_\text{soln}=10^{-7}$ \\
Truncation tolerance for convection matrix & $\epsilon_\text{conv}=10^{-3}$ \\
\hline
\end{tabular}
\end{table}

We demonstrate the savings obtained from the inexact solves. \Cref{table:inexact_picard} shows the performance of Picard's method if different stopping tolerances are used in \cref{eq:inexact_solve}. With a larger $tol_\text{gmres}$, the number of Picard steps does not increase, while the total number of GMRES iterations and the associated computational costs are greatly reduced.
\begin{table}[hbtp]
\label{table:inexact_picard}
\caption{Performance of Picard's method with different values of GMRES stopping tolerance $tol_\text{gmres}$. Truncation tolerance $\epsilon_\text{gmres}=10^{-2}*tol_\text{gmres}$. $tol_\text{picard}=10^{-5}$. LSC preconditioner is used.}
\centering
\begin{tabular}{c|ccc}
\hline
$tol_\text{gmres}$ & $10^{-1}$ & $10^{-3}$ & $10^{-5}$ \\
$\epsilon_\text{gmres}$ & $10^{-3}$ & $10^{-5}$ & $10^{-7}$ \\ \hline
Number of Picard steps & 5 & 5 & 5 \\
Total number of GMRES iterations & 18 & 39 & 58 \\
Computational time (s) & $205.9$ & $555.7$ & $1168.2$ \\
\hline
\end{tabular}
\end{table}

\begin{figure}[hbtp]
    \label{fig:prec}
    \begin{center}
    \subfloat[]{\label{fig:prec_num_iter}\includegraphics[width=0.5\textwidth]{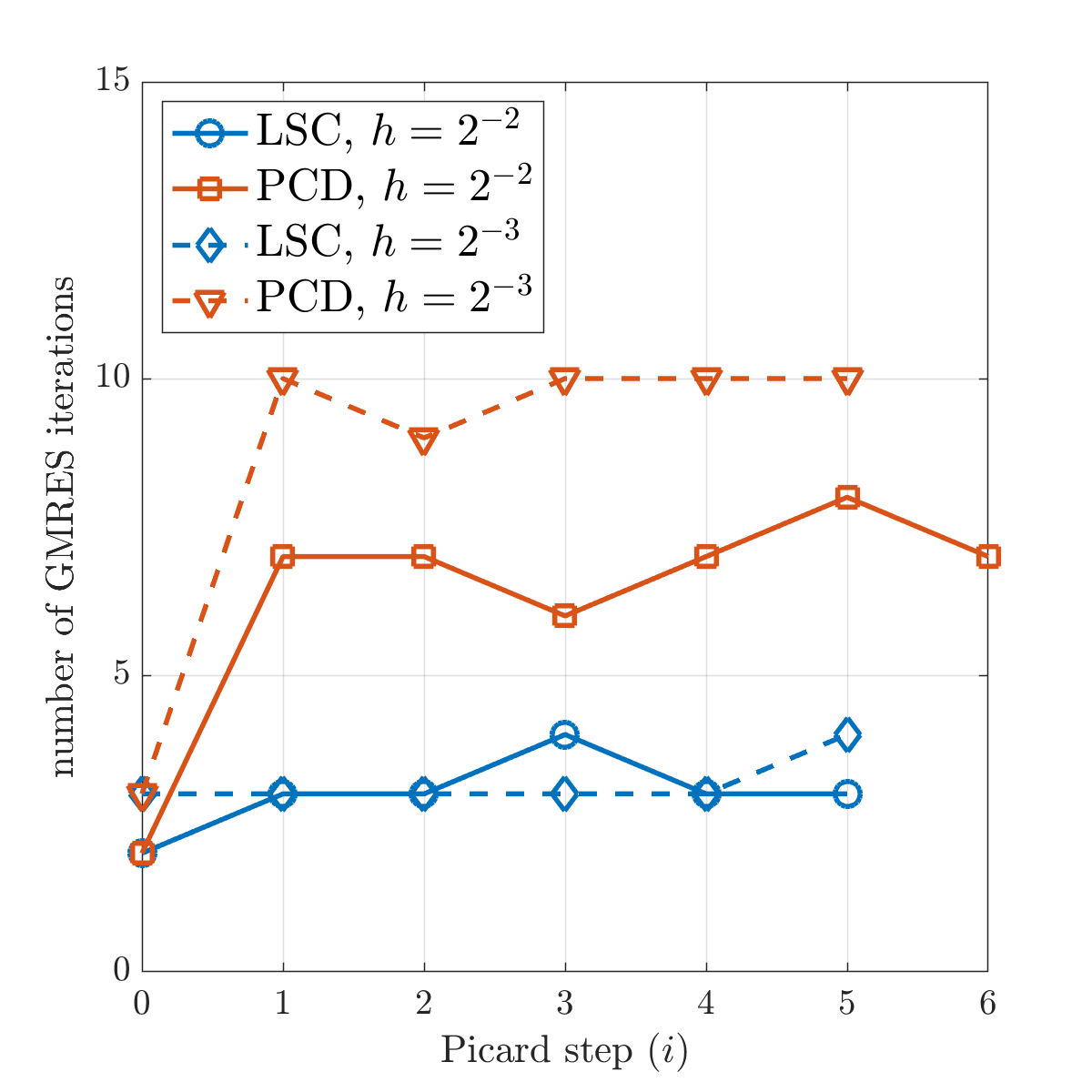}}
    \subfloat[]{\label{fig:prec_time}\includegraphics[width=0.5\textwidth]{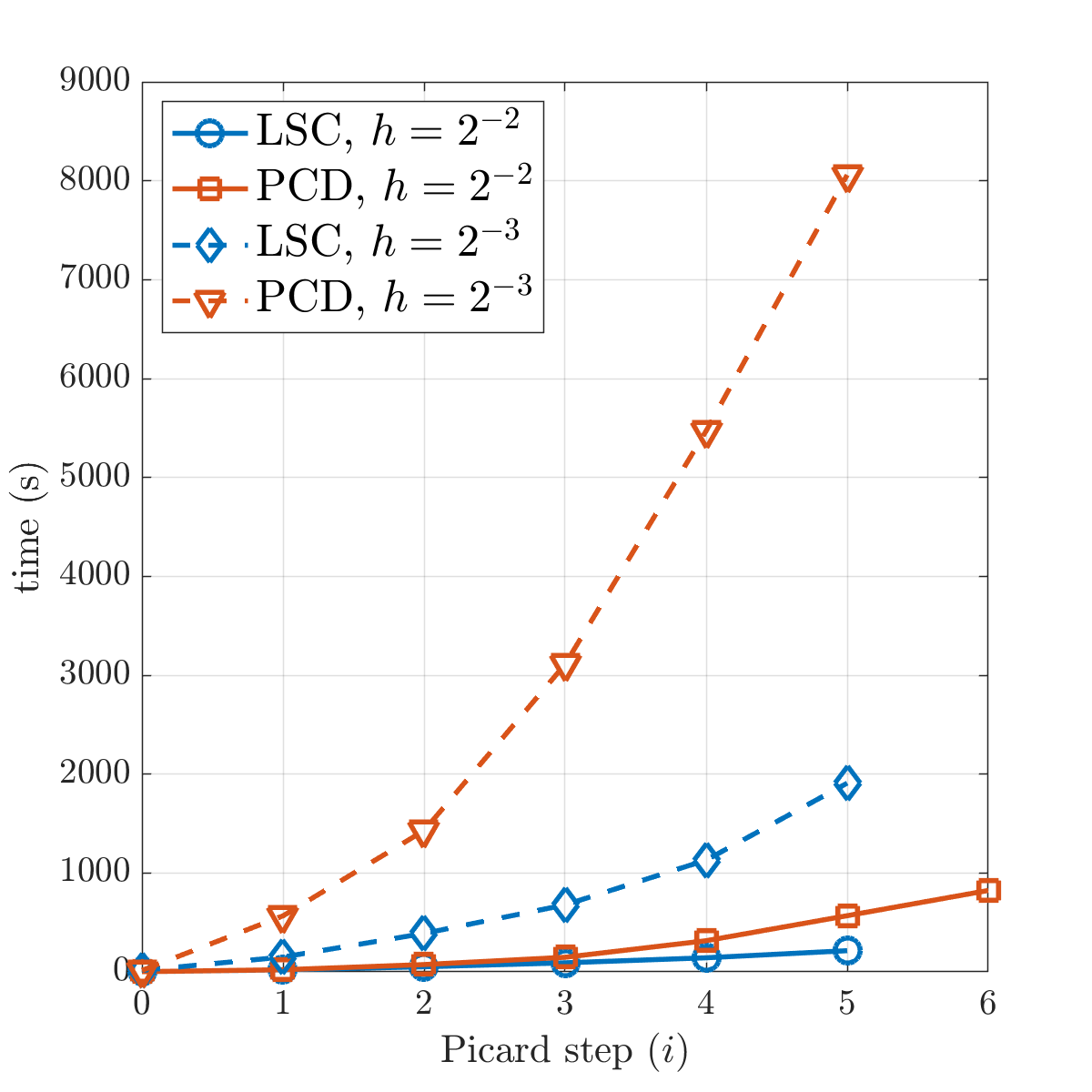}}
    \end{center}
    \caption{(a) Number of GMRES iterations at each Picard step. (b) Cumulative computational time after each Picard step. For $h=2^{-2}$ with PCD preconditioner, it takes 6 Picard steps for convergence.}
\end{figure}

We compare the two mean-based preconditioners discussed in \cref{sec:prec}. \Cref{fig:prec} shows the number of GMRES iterations required at each Picard step, and the associated computational costs when the two preconditioners are used. For two different mesh sizes, the PCD preconditioner results in larger numbers of GMRES iterations, and thus higher computational times, than the LSC preconditioner. It should also be noted that for both preconditioners, only a small number of GMRES iterations is needed for solving the linear system at each Picard step. This is partially due to the large stopping tolerance used in \cref{eq:inexact_solve}. The LSC preconditioner will be used for the numerical tests below.

\begin{figure}[hbtp]
    \label{fig:ranks_para}
    \begin{center}
    \subfloat[]{\label{fig:sigma}\includegraphics[width=0.5\textwidth]{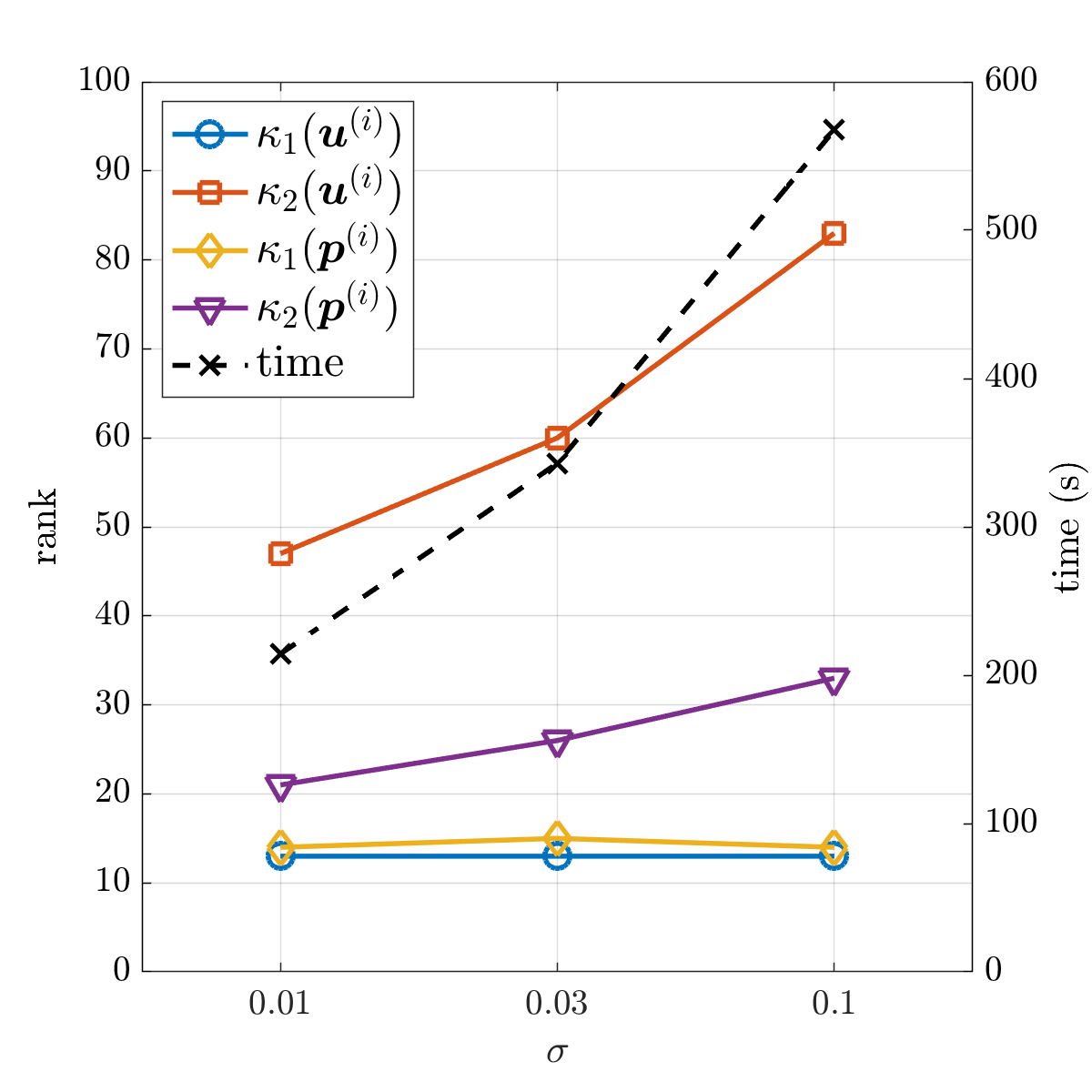}}
    \subfloat[]{\label{fig:nu0}\includegraphics[width=0.5\textwidth]{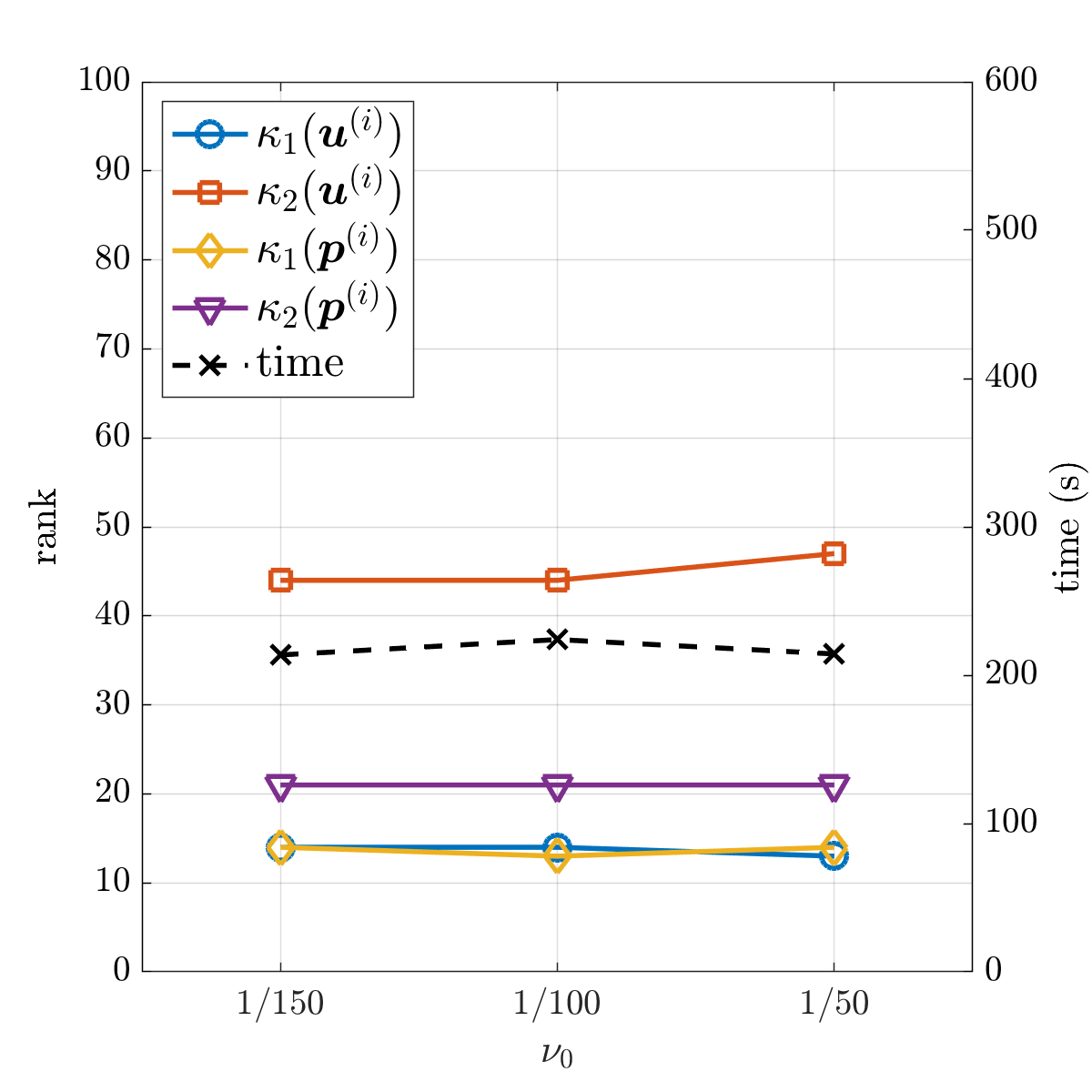}}
    \end{center}
    \caption{Tensor train ranks of solutions $\bm{u}^{(i)}$ and $\bm{p}^{(i)}$ at final Picard step and computational times to compute solutions, for different values of $\sigma$ and $\nu_0$.}
\end{figure}

In the following, we test the algorithm with several variants of the benchmark problem determined by various values of parameters associated with it. \Cref{fig:sigma} shows the solution ranks and computational times for three different values of $\sigma$. When $\sigma$ is smaller, the standard deviation is smaller, the discrete solution can be approximated by a tensor train with smaller ranks, and it is also less expensive to solve the nonlinear problem. On the other hand, even for $\sigma=0.1$, the low-rank solution takes much less storage than a full tensor. For example, the ranks of the approximate solution $\bm{u}^{(i)}$ are $\kappa_1=13$, $\kappa_2=83$. The ratio of storage requirements between such a tensor train and a full tensor is
\begin{equation}
    \frac{n_t\kappa_1+n_\xi\kappa_1\kappa_2+n_u\kappa_2}{n_t n_\xi n_u} = \frac{270748}{3829760} \approx 7.1\%.
\end{equation}
The same quantities are plotted in \cref{fig:nu0} for different values of the mean viscosity $\nu_0$. The ranks and computational times are not significantly affected by $\nu_0$.

\begin{figure}[hbtp]
    \label{fig:ranks_size}
    \begin{center}
    \subfloat[]{\label{fig:h}\includegraphics[width=0.5\textwidth]{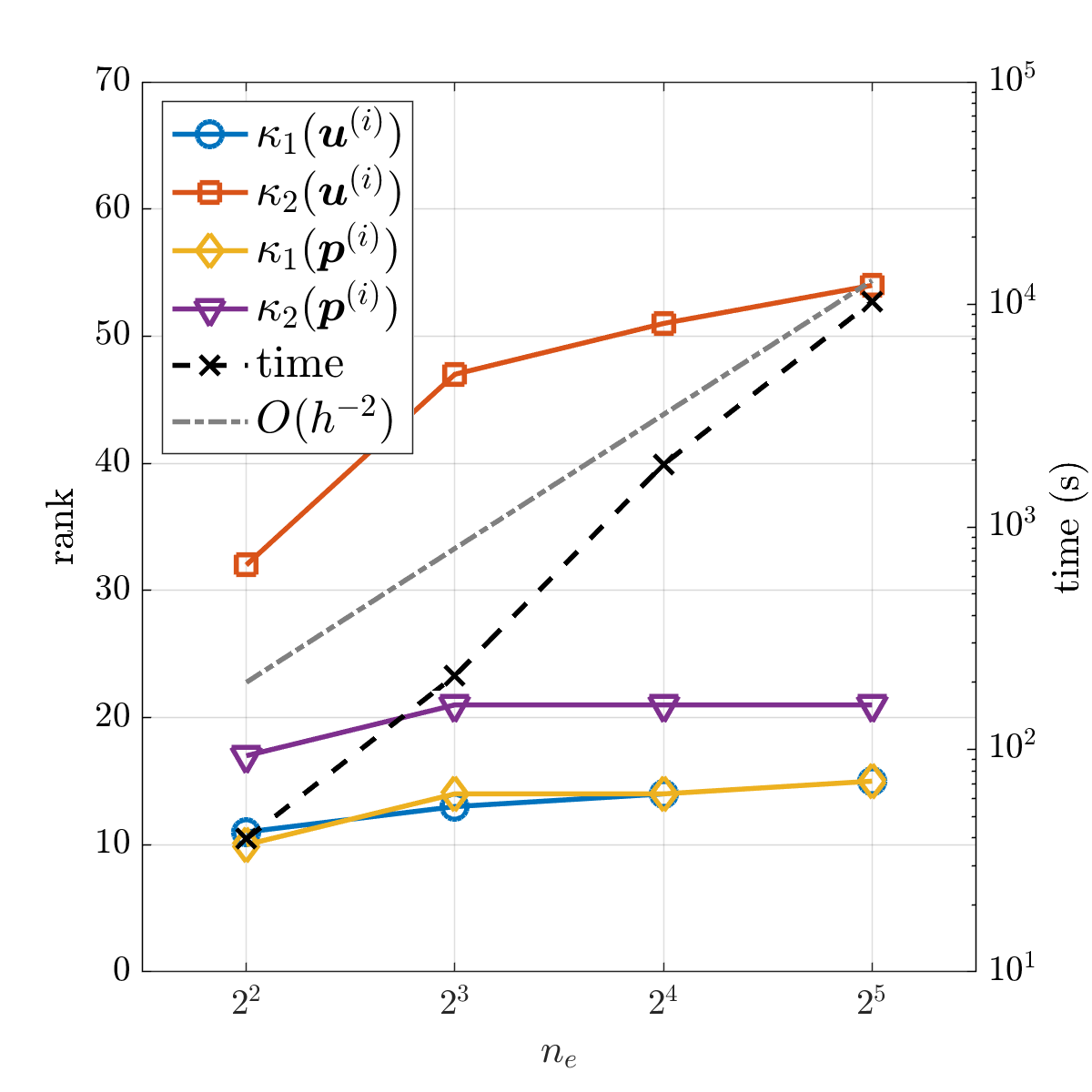}}
    \subfloat[]{\label{fig:tau}\includegraphics[width=0.5\textwidth]{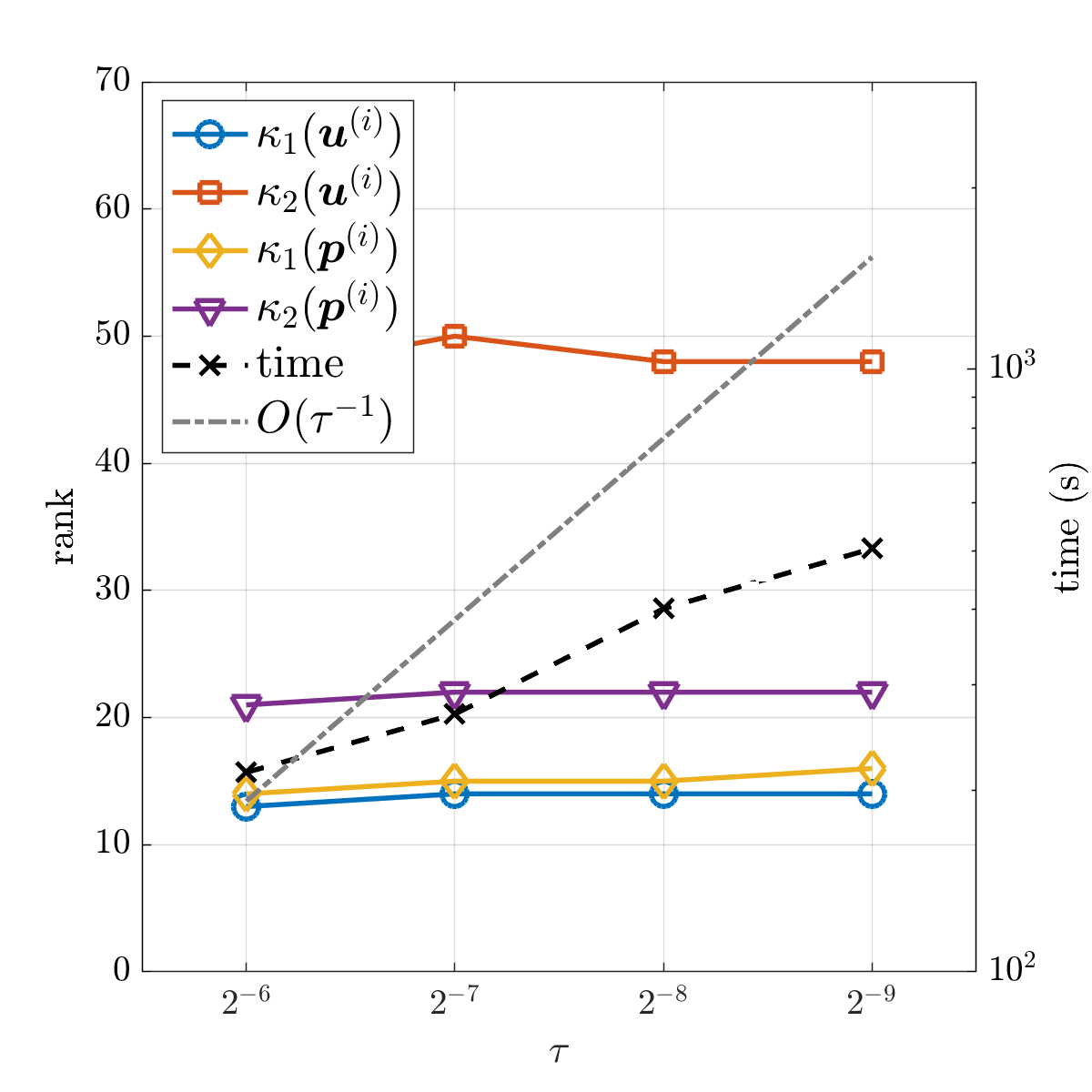}}
    \end{center}
    \caption{Tensor train ranks of solutions $\bm{u}^{(i)}$ and $\bm{p}^{(i)}$ at final Picard step and computational times to compute solutions, for different values of $h$ and $\tau$. In (a), $n_e=2/h$ is the number of elements in the vertical interval $[-1,1]$ of the domain $\mathcal{D}$.}
\end{figure}

Finally, the algorithm is applied to solve discrete problems with various mesh sizes $h$ or time step sizes $\tau$. It can be seen from \cref{fig:h} that there is only a slight increase in the solution ranks as the spatial mesh is refined. It is also shown in \cref{fig:h} that the computational time increases with an asymptotic rate $O(h^{-2})$ (note that a logarithmic scale is used in the figure). In other words, as the spatial mesh is refined, no extra computational burden is introduced except for the increased problem size. For different time step sizes $\tau$, the computational time increases much more slowly than $O(\tau^{-1})$ (see \cref{fig:tau}). This is due to the fact that in $\mathbb{F}+\mathbb{C}$, the matrices obtained from time discretization are very simple (e.g., $I_{n_t}$ and $C_{n_t}$), and thus an increase in $n_t$ does not make a significant impact on the computational costs.

%-------------------------------------------------------------
\section{Conclusions}
\label{sec:conclusion}

In this paper, we developed and studied efficient low-rank iterative methods for solving the time-dependent Navier--Stokes equations with a random viscosity. We considered an all-at-once formulation where the discrete solutions at all the time steps are solved together in a single system. To address the high storage and computational costs of this strategy, we used low-rank tensor approximations in a Newton--Krylov type algorithm. For the all-at-once system, we proposed two mean-based preconditioners using results from the deterministic problem. The computational costs were further reduced with inexact Picard method and approximate convection matrices. It was shown in the numerical experiments that the low-rank method is able to solve the nonlinear problem efficiently and the discrete solutions have small tensor ranks.

\bibliographystyle{plain}
\bibliography{stoch_ns}

\end{document}